\pdfoutput=1
\documentclass[pdflatex,sn-mathphys-ay]{sn-jnl}

\usepackage{graphicx}%
\usepackage{amsmath,amssymb,amsfonts}%
\usepackage{mathtools,bm}%
\usepackage{mathrsfs}%
\usepackage[title]{appendix}%
\usepackage{xcolor}%
\usepackage{textcomp}%
\usepackage{manyfoot}%
\usepackage{booktabs,array,tabularx}%
\usepackage{algorithm}%
\usepackage{algorithmicx}%
\usepackage[noend]{algpseudocode}%
\usepackage{enumitem}%
\setlist{itemsep=2pt,topsep=4pt}

\newcommand{\Nor}{\mathcal N}
\newcommand{\R}{\mathbb R}
\newcommand{\E}{\mathbb E}
\newcommand{\T}{^{\top}}
\newcommand{\cS}{\mathcal S}
\newcommand{\Hb}{\bar H}
\DeclareMathOperator{\range}{range}
\DeclareMathOperator{\diag}{diag}
\DeclareMathOperator{\Var}{Var}
\DeclareMathOperator{\orth}{orth}
\DeclareMathOperator{\tr}{tr}
\DeclareMathOperator{\gap}{gap}
\DeclareMathOperator{\IACT}{IACT}
\DeclareMathOperator{\polar}{polar}

\theoremstyle{thmstyleone}%
\newtheorem{proposition}{Proposition}%
\theoremstyle{thmstyletwo}%

\begin{document}

\title[Gradient-free MCMC exploiting Gaussian and non-Gaussian posterior structure]{Gradient-free MCMC exploiting low-dimensional Gaussian and non-Gaussian posterior structure}

\author[1]{\fnm{Hongqiao} \sur{Wang}}\email{Hongqiao.Wang@csu.edu.cn}
\author[2]{\fnm{Jiahua} \sur{Jiang}}\email{jiangjh@shanghaitech.edu.cn}
\author[3]{\fnm{Jing} \sur{Liu}}\email{liujingdx@ncwu.edu.cn}
\author*[4]{\fnm{Jinglai} \sur{Li}}\email{j.li.10@bham.ac.uk}

\affil[1]{\orgdiv{School of Mathematics and Statistics}, \orgname{Central South University}, \orgaddress{\city{Changsha}, \postcode{410083}, \country{China}}}
\affil[2]{\orgdiv{Institute of Mathematical Sciences}, \orgname{ShanghaiTech University}, \orgaddress{\city{Shanghai}, \country{China}}}
\affil[3]{\orgdiv{College of Water Resources}, \orgname{North China University of Water Resources and Electric Power}, \orgaddress{\city{Zhengzhou}, \country{China}}}
\affil*[4]{\orgdiv{School of Mathematics}, \orgname{University of Birmingham}, \orgaddress{\street{Edgbaston}, \city{Birmingham}, \postcode{B15 2TT}, \country{UK}}}

\abstract{Markov chain Monte Carlo (MCMC) for Bayesian inverse problems is expensive: the unknown is high-dimensional, every likelihood evaluation requires a forward solve, and gradients are often unavailable. The posterior typically differs from the prior in a few directions and is often approximately Gaussian in most of them but strongly non-Gaussian in a few. We propose the NQU sampler, a gradient-free MCMC method that splits the parameter space into non-Gaussian (N), approximately quadratic (Q) and uninformative (U) directions. The split is learned once, before sampling, from likelihood values, and each block is updated by a kernel suited to it; a collapsed update handles a global prior scale. Because the learned structure only shapes the proposals, every update leaves the exact posterior invariant, at the cost of one likelihood evaluation. We give a condition under which the blocks can be learned and show that, for likelihoods with exactly the learned structure, the spectral gap of the sampler does not depend on the dimension, whereas updating the scale given all coordinates can have a gap of order $1/d$. On the posterior axes, the NQU sampler, learning included, needed the fewest evaluations in four of five synthetic examples with $d=200$, among nine baselines including gradient-based ones; in those with non-Gaussian directions, gradient-free baselines needed at least ten times as many. On a methane ignition problem with 219 uncertain rate constants, it gave consistent estimates, whereas full-space robust adaptive Metropolis needed about 300 times as many evaluations and an ensemble sampler degenerated.}

\keywords{Markov chain Monte Carlo, gradient-free sampling, Bayesian inverse problems, likelihood-informed subspace, Hermite expansion, adaptive MCMC}

\maketitle

\section{Introduction}\label{sec:intro}

Inverse problems have important applications across science and engineering, including medical imaging, geophysics, subsurface flow and chemical kinetics \citep{tarantola2005}. The task is to infer an unknown $\theta$ from finitely many indirect and noisy observations $y=G(\theta)+\epsilon$, where evaluating the forward model $G$ typically means solving a partial or ordinary differential equation, and $\theta$ is, for example, a spatially varying coefficient, an initial condition or a set of physical parameters. Such problems are usually ill-posed: the data alone determine $\theta$ neither uniquely nor stably. The Bayesian approach places a prior on $\theta$, most often a Gaussian $\Nor(0,C_0)$, and takes the posterior distribution as the solution \citep{kaipio2005,stuart2010,dashti2017}. The prior regularises the problem, and the posterior quantifies the remaining uncertainty, for example through credible intervals for $\theta$ and for predictions made with it. With $\Phi(\theta)$ the negative log-likelihood, the posterior has density proportional to $\exp\{-\Phi(\theta)\}$ with respect to the prior; for nonlinear $G$ it has no closed form, and Markov chain Monte Carlo (MCMC), which is asymptotically exact, is the standard tool for sampling from it \citep{robert2004}. In inverse problems, however, MCMC is expensive. The main reason is the dimension. A spatially varying unknown must be discretised, so that $\theta\in\R^d$ with $d$ growing as the mesh is refined, and a chemical reaction mechanism can have hundreds of uncertain rate constants. As $d$ grows, random-walk proposals must take smaller steps, successive samples become more strongly correlated, and more steps are needed for a given number of effective samples \citep{roberts1997}. Moreover, each sample requires an evaluation of the forward model, which in many practical inverse problems is computationally demanding, so that samples come at a high cost.

Posteriors in inverse problems often have a structure that a sampler can exploit. Many practical inverse problems are ill-posed in the sense that the data do not provide information along every direction of the parameter space. In many of them the data are more informative than the prior in only a few directions, because the prior correlates the components of $\theta$, the observations are few and noisy, and the forward model is smoothing \citep{flath2011,buithanh2013,spantini2015,cui2014,zahm2022}. The posterior then differs from the prior essentially only in these directions. Many MCMC methods find a subspace containing these directions by various dimension-reduction techniques and use it to accelerate sampling. Most of them use gradients of the log-likelihood, and some also Hessians. Stochastic Newton MCMC \citep{martin2012} builds Gaussian proposals from low-rank approximations of the prior-preconditioned Hessian at the current state. The likelihood-informed subspace (LIS) of \citet{cui2014} is spanned by the leading eigenvectors of the prior-preconditioned Gauss--Newton Hessian averaged over the posterior; MCMC is run in this subspace, and the prior is used in its complement. Dimension-independent likelihood-informed (DILI) MCMC \citep{cui2016} samples the exact posterior with Langevin proposals that treat the LIS and its complement differently and remain efficient as the discretisation is refined. Geometric MCMC methods \citep{beskos2017,lan2019} use Langevin and Hamiltonian proposals informed by gradients and curvature in such a subspace, and preconditioned Crank--Nicolson (pCN) proposals \citep{cotter2013} in its complement. Active subspaces \citep{constantine2014} are spanned by the leading eigenvectors of the outer product of the log-likelihood gradient averaged over the prior, and AS-MCMC \citep{constantine2016} runs MCMC on the active variables.

In many practical problems, however, gradients may not be easy to obtain. The forward model is often a legacy or commercial code that can only be run as a black box. Even when the source code is available, computing gradients efficiently requires an adjoint solver, which takes considerable effort to derive and implement for each model, and codes with adaptive time stepping or switching logic may not be differentiable. To our knowledge, only a few dimension-reduction MCMC methods are gradient-free. Two of them use a subspace fixed in advance, spanned by the leading Karhunen--Lo\`eve modes of the prior, with pCN in its complement: \citet{zhou2017} run adaptive Metropolis in this subspace, and the functional ensemble sampler of \citet{coullon2021} moves an ensemble of walkers in it by affine-invariant steps. SAFES \citep{dunlop2022} instead adapts the subspace to the current ensemble. \citet{beskos2018} replace the Hessians in the DILI proposals by the empirical covariance of the particles in a multilevel sequential Monte Carlo method. Active subspaces can also be estimated without gradients, from a global quadratic model fitted by least squares \citep{constantine2017}. Such a model, however, has $(d+1)(d+2)/2$ coefficients, so the number of likelihood evaluations needed to fit it grows quadratically with $d$.

Another issue concerns the Gaussianity of the posterior. The efficiency of the MCMC methods above, gradient-based or gradient-free, relies on the posterior in the subspace being close to Gaussian or only mildly non-Gaussian. They identify the subspace from second-order information, such as averaged Hessians or covariances, which does not distinguish strongly non-Gaussian directions from approximately Gaussian ones, and they use a single type of proposal, mostly Gaussian or Langevin-type, for the whole subspace. The assumption of an approximately Gaussian posterior is justified to some extent. The posterior is exactly Gaussian when $G$ is linear and the noise is Gaussian, and for nonlinear $G$ it approaches a Gaussian as the noise level decreases, under conditions such as a unique mode \citep{lu2017,schillings2020}; this is also the basis of the common use of Laplace approximations in large-scale inverse problems. In some problems, however, although the posterior is approximately Gaussian in most directions, it can be strongly non-Gaussian in a few. Symmetries or non-identifiability produce several modes, one-sided constraints or thresholds produce skewness, heavy-tailed noise produces spikes, and hyperparameters such as the scale of the prior produce funnel-shaped posteriors \citep{agapiou2014}. These few directions can then make the whole chain mix slowly: its steps must be small enough for a spike or the neck of a funnel, and it rarely moves between modes.

In this paper we propose a gradient-free, structure-exploiting MCMC method that handles such strongly non-Gaussian directions. Besides a Gaussian prior, which can be relaxed (Section~\ref{sec:tau}), it assumes that the likelihood is informative only in a low-dimensional subspace, that only a few directions in this subspace are strongly non-Gaussian, and that the others are approximately Gaussian given these (Section~\ref{sec:setup}); these assumptions affect efficiency only, not correctness. Accordingly, the method splits the parameter space into three blocks: within the data-informed subspace, the non-Gaussian (N) directions and the approximately quadratic (Q) directions, in which the negative log-likelihood is nearly quadratic given the N block; and outside it, the uninformative (U) directions. We refer to it as the NQU sampler, after these blocks.

The blocks are learned once, before sampling, from likelihood values at independent design points, with a low-rank quadratic model for the data-informed subspace and third- and fourth-order moments for its non-Gaussian directions, in the spirit of principal Hessian directions \citep{li1992} and non-Gaussian component analysis \citep{blanchard2006}. The learned structure is frozen and used only to construct proposals, so every update leaves the exact posterior invariant. Sampling updates each block with a kernel suited to it: a random walk for the few non-Gaussian directions, an adaptive Gaussian proposal for the approximately quadratic ones and prior-based proposals for the uninformative ones. Each update costs one evaluation of $\Phi$, or \emph{call}. We give a condition under which the blocks can be learned from likelihood values and show that, for likelihoods with exactly the learned structure, the spectral gap of the sampler does not depend on $d$, and that the collapsed scale update avoids the gap of order $1/d$ of updating the scale given all coordinates. In experiments on five synthetic examples, one of them a backward heat conduction problem, and on a methane ignition problem, counting all calls including those for learning, the NQU sampler needed fewer calls than gradient-based baselines in most cases and, in the examples with non-Gaussian directions, at least an order of magnitude fewer than gradient-free ones. A Python implementation of the NQU sampler is available at \url{https://github.com/lijinglai/nqu-sampler}.

\section{The NQU sampler}\label{sec:method}

\subsection{Setting and assumptions}\label{sec:setup}

Let the prior be $\Nor(0,C_0)$ and let $x=C_0^{-1/2}\theta$, so that $x\sim\Nor(0,I_d)$ a priori. Writing $\Phi(x)$ for the negative log-likelihood, the target is
\[
\pi(x)\ \propto\ \exp\big\{-\tfrac12|x|^2-\Phi(x)\big\}.
\]
Non-Gaussian priors are treated in Section~\ref{sec:tau}. Let $\cS^*$ be the smallest subspace with $\Phi(x)=\Phi(P_{\cS^*}x)$, the \emph{informed subspace}. A direction $v$ in $\cS^*$ along which the data are less informative than the prior is called \emph{weakly informed}; for a likelihood term $\tfrac12\lambda(v\T x-m)^2$ this means $\lambda<1$, the prior precision. A \emph{structure} consists of matrices $W_N\in\R^{d\times r}$ and $W_Q\in\R^{d\times k}$ with orthonormal columns and $W_N\T W_Q=0$, and the coordinates
\[
a=W_N\T x\in\R^r,\qquad c=W_Q\T x\in\R^k,\qquad \rho=P_Ux,\quad P_U=I_d-W_NW_N\T-W_QW_Q\T ,
\]
called the N, Q and U blocks. They are independent under the prior and $|x|^2=|a|^2+|c|^2+|\rho|^2$. The N block ($r$ small) holds the directions in which the likelihood is markedly non-Gaussian; the Q block holds informed directions in which, given $a$, the likelihood is approximately quadratic; the U block holds the approximately uninformative directions. Unless stated otherwise, $\E$ refers to $x\sim\Nor(0,I_d)$.

Under $\Nor(0,I_d)$, $\Phi\in L^2$ has the Hermite (Wiener chaos) expansion $\Phi=\sum_{k\ge0}\Phi_k$, where $\Phi_k$ is a combination of $\mathrm{He}_\alpha(x)=\prod_i\mathrm{He}_{\alpha_i}(x_i)$ with $|\alpha|=k$. For a unit vector $v$ define the $k$th-order \emph{energy along $v$}
\[
E_k(v)=\|\Phi_k\|^2-\big\|\Phi_k^{(v^\perp)}\big\|^2 ,
\]
where $\Phi_k^{(v^\perp)}$ collects the terms with $\alpha_v=0$ in an orthonormal basis containing $v$; mixed terms count for each direction they involve. Then $E_{\rm tot}(v)=\sum_{k\ge1}E_k(v)=\E[\Var(\Phi\mid P_{v^\perp}x)]$, and we write $E_{\ge3}(v)=\sum_{k\ge3}E_k(v)$. Let
\[
g=\E[\Phi(x)\,x],\qquad \Hb=\E\big[\Phi(x)(xx\T-I_d)\big],
\]
which for smooth $\Phi$ are the average gradient and Hessian by Stein's identity \citep{stein1981}.

The following assumptions affect efficiency only, not correctness (Section~\ref{sec:props}).
\begin{description}[leftmargin=1.5em,style=nextline]
\item[A1 (small N block)] The number $r$ of non-Gaussian directions is at most a preset $r_{\max}$ (here 5).
\item[A2 (low-order dependence)] For every unit vector $v$ in the N block and every unit vector $v$ in the Q block, $E_{\ge3}(v)\le E_1(v)+E_2(v)$.
\item[A3 (structure of Q and U)] On the posterior scale, given $a$, $\Phi(x)\approx\phi_0(a)+(\ell+Ba)\T c+\tfrac12c\T Hc$, and $\Phi$ depends on $x$ approximately only through $(a,c)$.
\end{description}
If A3 holds exactly, Q directions satisfy A2 automatically, so A2 is mainly a condition on the N block.

\subsection{Learning the blocks}\label{sec:learn}

Learning uses $n_1+n_2$ calls, here $n_1=30\,000$ and $n_2=20\,000$, in three steps. The designs of Steps~1 and~2 are independent, so zero-mean noise in $\Phi$ changes the variance but not the population targets.

Step~1 (quadratic sensing)\label{sec:step1} draws $x_i\sim\Nor(0,I_d)$, $i=1,\dots,n_1$, evaluates $y_i=\Phi(x_i)$ and, holding out a validation set, fits the signed low-rank quadratic model
\[
q(x)=\kappa_0+\hat g\T x+\tfrac12x\T\hat Hx,\qquad \hat H=LDL\T=\sum_{j=1}^sd_j\,l_jl_j\T,\quad d_j\in\{-1,+1\},
\]
where $(\kappa_0,\hat g)$ is fitted by least squares on all $d$ coordinates and $L$ minimises the squared error with the linear part profiled out. The rank grows one at a time. For the current residuals $e_i$, the principal Hessian direction (PHD) matrix $\hat M=\tfrac12\big[n^{-1}\sum_ie_ix_ix_i\T-\bar e\,I_d\big]$ estimates half the average residual Hessian; for its eigenpair $(\lambda,v)$ of largest $|\lambda|$ we add the column $l=\sqrt{2|\lambda|}\,v$ with sign $\operatorname{sign}\lambda$ and refine all columns by scaled gradient descent \citep{tong2021} with fixed signs. We stop, discarding the new column, when it reduces the validation error by less than $\epsilon_0=2\%$. The matrix $\hat H$ is not restricted to be positive semidefinite, because negative average curvature is informative: for $\Phi=5e^{-z^2}$, $z=v\T x$, we have $g=0$ but $v\T\Hb v\approx-1.92$. The output $\cS$ is the span of the eigenvectors of $\hat g\hat g\T+\hat H^2$ with eigenvalues above $\delta^2$, $\delta=0.1$; let $m=\dim\cS$ and $F_{\cS}\in\R^{d\times m}$ an orthonormal basis. Directions outside $\cS$ form the U block; Section~\ref{sec:step1pop} analyses what this step recovers.

Step~2\label{sec:step2} selects the non-Gaussian block from a second design. Step~1 gives on $\cS$ the reference precision $A_{\cS}=I_m+F_{\cS}\T\hat HF_{\cS}$, with its eigenvalues bounded below (Table~\ref{tab:impl}), and mean $\mu_{\cS}=-A_{\cS}^{-1}F_{\cS}\T\hat g$. Without the bound this is the prior times $e^{-q}$ on $\cS$; along directions where the bound is active, the reference is wider than the prior. With $w=A_{\cS}^{1/2}(F_{\cS}\T x-\mu_{\cS})$, draw
\[
x_i=F_{\cS}\big(\mu_{\cS}+A_{\cS}^{-1/2}w_i\big)+P_{\cS^\perp}\xi_i,\qquad w_i\sim\Nor(0,I_m),\quad \xi_i\sim\Nor(0,I_d),
\]
for $i=1,\dots,n_2$,
and evaluate $\Phi(x_i)$. The centred residuals $e_i=\Phi(x_i)-q(x_i)$ have the same third- and fourth-order Hermite coefficients in $w$ as $\Phi$, because $q$ is quadratic in $(w,P_{\cS^\perp}\xi)$; when $q$ fits well, subtracting it reduces the variance of their estimates.

To detect non-Gaussian directions, split the design into halves $h=1,2$ and estimate $\mathcal T_3^{(h)}=\widehat\E_h[e\,\mathrm{He}_3(w)]\in\R^{m\times m^2}$ and $\mathcal T_4^{(h)}=\widehat\E_h[e\,\mathrm{He}_4(w)]\in\R^{m\times m^3}$, unfolded along the first index. The cross-fitted matrix
\[
M_{3,4}=\frac{\mathcal T_3^{(1)}\mathcal T_3^{(2)\top}+\mathcal T_3^{(2)}\mathcal T_3^{(1)\top}}{2\cdot3!}+\frac{\mathcal T_4^{(1)}\mathcal T_4^{(2)\top}+\mathcal T_4^{(2)}\mathcal T_4^{(1)\top}}{2\cdot4!}
\]
has no positive noise bias: $\E M_{3,4}=\mathcal T_3\mathcal T_3\T/3!+\mathcal T_4\mathcal T_4\T/4!$ with $\mathcal T_k=\E[\Phi\,\mathrm{He}_k(w)]$. In a basis containing a unit vector $u$, $u\T(\E M_{3,4})u$ counts the energy of each $\mathrm{He}_\alpha$ with weight $\alpha_u/|\alpha|$, so the leading eigenvector carries the most third- and fourth-order energy; the fourth order is needed because $\mathcal T_3\equiv0$ for even nonlinearities.

Directions are selected sequentially. With $B$ an orthonormal basis of the directions not yet selected ($B=I_m$ initially), the statistic is $\lambda_{\max}(B\T M_{3,4}B)$. Its null distribution is approximated by random sign flips $e_i\mapsto s_ie_i$, $s_i=\pm1$, which keep the pairing of $|e_i|$ with $w_i$ but break the dependence of the sign of $e_i$ on $w_i$. If the statistic exceeds a high quantile of the flips (Table~\ref{tab:impl}), the leading eigenvector of $B\T M_{3,4}B$, mapped back to $\R^m$, is selected and the test is repeated, for at most $r_{\max}=5$ directions. The selected directions form $E\in\R^{m\times r}$.

To map the selection back to $x$, let $W_N=\orth(F_{\cS}A_{\cS}^{-1/2}E)$ and let $W_Q$ be an orthonormal basis of the orthogonal complement of $W_N$ in $\cS$, $k=m-r$. Then $c=W_Q\T x$ is an invertible affine function of the unselected coordinates $E_\perp\T w$, and moving along $W_N$ changes only the selected coordinates of $w$: the columns of $W_N$ are the displacement directions of the selected non-Gaussian functionals under the reference Gaussian. If N and Q are coupled ($B\ne0$ in A3), these directions contain components of the coupled Q directions, and the learned coordinates are only approximately exact (Appendix~\ref{app:backmap}).

Since (N) is a coordinatewise random walk, its efficiency depends on the basis of $\range(W_N)$: a coordinate mixing a bimodal and a spiky direction takes steps limited by the spike. Using the same data (no calls), let $z=E\T w$ and $\mathcal C_k^{(h)}=\widehat\E_h[e\,\mathrm{He}_k(z)]$, $k=3,4$, and maximise over orthogonal $R\in O(r)$, with columns $R_j$, the cross-fitted contrast
\[
\mathcal J(R)=\sum_{j=1}^r\Big[\frac{\mathcal C_3^{(1)}(R_j^{\otimes3})\,\mathcal C_3^{(2)}(R_j^{\otimes3})}{3!}+\frac{\mathcal C_4^{(1)}(R_j^{\otimes4})\,\mathcal C_4^{(2)}(R_j^{\otimes4})}{4!}\Big]
\]
by Jacobi rotations. $\mathcal J$ is the energy of the pure third- and fourth-order terms of the rotated coordinates; as in independent component analysis \citep{cardoso1993,comon1994}, maximising it makes the two tensors as diagonal as possible. With $\gamma_j=F_{\cS}A_{\cS}^{1/2}(ER)_j$ the corresponding functional of $x$, we replace $W_N$ by $W_N\polar(W_N\T\Gamma)$, $\Gamma=[\gamma_1/|\gamma_1|,\dots,\gamma_r/|\gamma_r|]$, the orthonormal basis of $\range(W_N)$ closest to the $\gamma_j$.

Step~3\label{sec:step3} builds the Q model, the part of the Step~1 model on $(a,c)$: $\ell=W_Q\T\hat g$, $B=W_Q\T\hat HW_N$, $H=W_Q\T\hat HW_Q$, with the eigenvalues of $I_k+H$ restricted to $[0.5,10^6]$. Let $V=(I_k+H)^{-1}$ and $\mu(a)=-V(\ell+Ba)$, and fix a square root $V^{1/2}$. Under A3, given $(a,\rho)$, $c$ is approximately $\Nor(\mu(a),V)$. This step needs no calls; $W_N,W_Q,\ell,B,H$ are then frozen.

\subsection{Sampling with the learned structure}\label{sec:sample}

Write $x=W_N a+W_Q c+\rho$ with $c=\mu(a)+V^{1/2}u$. The map $(a,u,\rho)\mapsto x$ is affine with constant Jacobian, and under exact A3 the conditional distribution of $u$ given $(a,\rho)$ is $\Nor(0,I_k)$. Each update evaluates $\Phi$ once, at the proposal; $\Phi$ at the current state is cached.

The (N) update chooses $j\in\{1,\dots,r\}$ uniformly and sets $a_j'=a_j+\sigma_j\xi$, $\xi\sim\Nor(0,1)$, keeping the other coordinates, $u$ and $\rho$ fixed:
\[
c'=\mu(a')+V^{1/2}u=c-VB(a'-a),\qquad x'=W_N a'+W_Q c'+\rho .
\]
It accepts with probability $\alpha_N=\min\{1,\exp(-\tfrac12|a'|^2-\tfrac12|c'|^2-\Phi(x')+\tfrac12|a|^2+\tfrac12|c|^2+\Phi(x))\}$. Because the Q block follows its learned conditional mean, the N block need not wait for it; under exact A3, (N) is a Metropolis--Hastings step on the marginal posterior of $a$ (Proposition~\ref{prop:exact}).

The (Q) update keeps $a$ and $\rho$ fixed, proposes by RAM \citep{vihola2012}, $c'=c+S_{n-1}\xi$ with $\xi\sim\Nor(0,I_k)$, accepts with probability $\alpha_n=\min\{1,\exp(-\tfrac12|c'|^2+\tfrac12|c|^2-\Phi(x')+\Phi(x))\}$, and updates
\[
S_nS_n\T=S_{n-1}\Big(I_k+\eta_n(\alpha_n-\alpha^*)\frac{\xi\xi\T}{|\xi|^2}\Big)S_{n-1}\T,\qquad \eta_n=\min\{1,kn^{-2/3}\},\quad\alpha^*=0.234,
\]
with $S_n$ lower triangular, $n$ the number of (Q) updates so far and $S_0S_0\T=(2.38^2/k)V$ \citep{roberts1997}. RAM adapts throughout sampling. The Q model enters only through $S_0$ and the shift in (N); its errors do not enter the acceptance of (Q).

The (U) update keeps $a$ and $c$ fixed. Draw $\xi\sim\Nor(0,I_d)$; with probability $\tfrac12$ set $\rho'=P_U\xi$ (refresh from the conditional prior), otherwise $\rho'=\sqrt{1-b_U^2}\,\rho+b_UP_U\xi$ (pCN step), and accept with probability $\alpha_U=\min\{1,\exp(\Phi(x)-\Phi(x'))\}$. Both proposals are reversible with respect to $\Nor(0,P_U)$, so uninformative directions do not lower the acceptance rate, and an informative direction left in U can still cross a barrier through an accepted refresh.

Updates are chosen with fixed probabilities 0.4, 0.3, 0.3 for (N), (Q), (U) (0.5 and 0.5 for (Q), (U) if $r=0$). The chain starts at $x_0=F_{\cS}\mu_{\cS}+P_{\cS^\perp}\xi$. The first 10\% of the sampling calls are a discarded warm-up in which the step sizes $\sigma_j$ and the pCN parameter $b_U$ are tuned towards target acceptance rates and then frozen (Table~\ref{tab:impl}). Algorithm~\ref{alg:main} summarises the method.

\begin{algorithm}[tbp]
\caption{NQU sampler}\label{alg:main}
\begin{algorithmic}[1]
\Require black box $\Phi$; learning calls $n_1,n_2$; sampling calls $T$, warm-up $T_w=\lfloor0.1T\rfloor$
\State \textbf{Step 1} ($n_1$ calls): $x_i\sim\Nor(0,I_d)$; fit $(\kappa_0,\hat g,\hat H)$; $\cS\gets$ leading eigenspace of $\hat g\hat g\T+\hat H^2$
\State \textbf{Step 2} ($n_2$ calls): design on $\Nor(\mu_{\cS},A_{\cS}^{-1})$; select $E$ by $M_{3,4}$ and sign flips; $W_N,W_Q$; if $r\ge2$, maximise $\mathcal J(R)$ and set $W_N\gets W_N\polar(W_N\T\Gamma)$
\State \textbf{Step 3}: $\ell,B,H$; $V$, $\mu(\cdot)$; $S_0\gets\mathrm{chol}\big((2.38^2/k)V\big)$
\State $x\gets F_{\cS}\mu_{\cS}+P_{\cS^\perp}\xi$; evaluate and cache $\Phi(x)$; $\sigma_j\gets1$, $b_U\gets1$
\For{$t=1,\dots,T$}
  \State choose $M\in\{N,Q,U\}$ with probabilities $(0.4,0.3,0.3)$; propose $x'$ as in Section~\ref{sec:sample}
  \State evaluate $\Phi(x')$ \Comment{one call}
  \State accept with probability $\alpha_M$; if $M=Q$, update $S_n$
  \State \textbf{if} $t\le T_w$ \textbf{then} tune $\sigma_j$ and $b_U$ \textbf{else} record $x$
\EndFor
\end{algorithmic}
\end{algorithm}

\subsection{Extension to non-Gaussian priors}\label{sec:tau}

A prior that is the image of $\Nor(0,I_d)$ under a known map $T$, for example a Laplace prior obtained by componentwise quantile transformations, is handled by sampling the Gaussian variable with $\Phi\circ T$ in place of $\Phi$ \citep[Remark~2.1]{dunlop2022}; the assumptions of Section~\ref{sec:setup} then concern $\Phi\circ T$, and a nonlinear $T$ can spread the dependence of the likelihood over many coordinates.

For a prior that is Gaussian given an unknown global scale, let $\theta\mid\tau\sim\Nor(0,\tau C_0)$ and $\tau\sim p(\tau)$, and normalise $\E\tau=1$. Then $x\mid\tau\sim\Nor(0,\tau I_d)$ and the target is
\[
\pi(x,\tau)\ \propto\ p(\tau)\,\tau^{-d/2}\exp\Big\{-\frac{|x|^2}{2\tau}-\Phi(x)\Big\}.
\]
Learning is unchanged, with design $\Nor(0,C_0)$, except that the negative eigenvalues of $H$ are set to zero, giving $H_+$, so that $I_k/\tau+H_+$ is positive definite for every $\tau$; the removed part is absorbed in the residual $R$ below.

Given $\tau$, the blocks $a,c,\rho$ are independent a priori with variances $\tau$. Let $b(a)=\ell+Ba$, $\Phi_Q(a,c)=b(a)\T c+\tfrac12c\T H_+c$, $R(x)=\Phi(x)-\Phi_Q(a,c)$,
\[
V_\tau=\big(I_k/\tau+H_+\big)^{-1},\qquad \mu_\tau(a)=-V_\tau b(a),\qquad u=V_\tau^{-1/2}\big(c-\mu_\tau(a)\big)
\]
with the symmetric square root, and
\[
Z_Q(a;\tau)=\int\Nor(c;0,\tau I_k)\,e^{-\Phi_Q(a,c)}\,dc=\Big(\frac{|V_\tau|}{\tau^k}\Big)^{1/2}\exp\big\{\tfrac12\,b(a)\T V_\tau b(a)\big\}.
\]
The updates (N), (Q), (U) are applied at fixed $\tau$ with the prior terms divided by $\tau$: (N) shifts $c'=c-V_\tau B(a'-a)$; (Q) runs RAM in $u$, $c'=c+V_\tau^{1/2}S_{n-1}\xi$ with $S_0=(2.38/\sqrt k)I_k$, which at fixed $\tau$ is equivalent to RAM in $c$ started from $(2.38^2/k)V_\tau$; and (U) uses $\sqrt\tau\,P_U\xi$ in place of $P_U\xi$.

The (H) update fixes $a$ and has three steps.
(1) Update $\eta=\log\tau$ by slice sampling with stepping out and shrinkage \citep{neal2003} (Table~\ref{tab:impl}) on
\[
\hat p(\eta\mid a)\ \propto\ p(e^\eta)\,e^\eta\,\Nor\big(a;0,e^\eta I_r\big)\,Z_Q\big(a;e^\eta\big),
\]
the conditional density of $\eta$ given $a$ when $\Phi$ is replaced by $\Phi_Q$ and $c,\rho$ are integrated out; this needs no call.
(2) Draw $u'\sim\Nor(0,I_k)$ and $\xi\sim\Nor(0,I_d)$, and set $c'=\mu_{\tau'}(a)+V_{\tau'}^{1/2}u'$, $\rho'=\sqrt{\tau'}\,P_U\xi$ and $x'=W_N a+W_Q c'+\rho'$.
(3) Accept $(x',\tau')$ with probability $\alpha_H=\min\{1,e^{R(x)-R(x')}\}$, one call.
Steps (1)--(2) form a proposal that is reversible with respect to the (possibly unnormalised; Appendix~\ref{app:phat}) distribution obtained by replacing $\Phi$ with $\Phi_Q$, so only the ratio $e^{-R}$ remains (Appendix~\ref{app:H}). Because (H) conditions only on $a$, $\tau$ need not wait for $c$ and $\rho$; if $\tau$ were drawn given all of $x$, its conditional density, proportional to $p(\tau)\tau^{-d/2}e^{-|x|^2/(2\tau)}$, would have relative width about $\sqrt{2/d}$. The selection probabilities are 0.36, 0.27, 0.27, 0.10 for (N), (Q), (U), (H) (0.45, 0.45, 0.10 for (Q), (U), (H) if $r=0$), and $\tau_0=1$.

\section{Theory}\label{sec:theory}

\subsection{Learnability}\label{sec:step1pop}

With $g$ and $\Hb$ of Section~\ref{sec:setup}, the first two energies are
\begin{equation}\label{eq:E12}
E_1(v)=(v\T g)^2,\qquad E_2(v)=|\Hb v|^2-\tfrac12(v\T\Hb v)^2 ,
\end{equation}
and the matrix $K=gg\T+\Hb^2$ captures them for every direction.

\begin{proposition}\label{prop:K}
Let $\Phi\in L^2(\Nor(0,I_d))$.
\begin{enumerate}[label=(\roman*)]
\item For every unit vector $v$, $\tfrac12\,v\T Kv\le E_1(v)+E_2(v)\le v\T Kv$.
\item If $v$ satisfies the inequality in A2, then $v\T Kv\ge\tfrac12E_{\rm tot}(v)$.
\item $\range(K)\subseteq\cS^*$, with equality if and only if no direction in $\cS^*$ has $E_1=E_2=0$; this holds if $\Hb$ restricted to $\cS^*$ is nonsingular.
\end{enumerate}
\end{proposition}

The proof is in Appendix~\ref{app:K}. Thus $K$ reflects most of the dependence of directions satisfying A2; if $E_{\ge3}(v)=t\,(E_1(v)+E_2(v))$, the first two orders carry only a fraction $1/(1+t)$ of it. That $K$ reflects a direction does not mean that the first learning step finds it (Proposition~\ref{prop:step1}).

To see what Step~1 learns, let $\Hb=\sum_j\mu_ju_ju_j\T$ with $|\mu_1|\ge|\mu_2|\ge\cdots$ ($\mu_j=0$ for $j>\operatorname{rank}\Hb$), and let $\Hb_s=\sum_{j\le s}\mu_ju_ju_j\T$ be its best rank-$s$ approximation. Allow observations $y=\Phi(x)+\xi$ with $\E[\xi\mid x]=0$ and $\E\xi^2=\sigma_\xi^2$, let $E^{\rm tot}_{\ge3}=\sum_{k\ge3}\|\Phi_k\|^2$ and
\[
\mathrm{MSE}_s=\tfrac12\sum_{j>s}\mu_j^2+E^{\rm tot}_{\ge3}+\sigma_\xi^2 .
\]

\begin{proposition}\label{prop:step1}
Let $\Phi\in L^2(\Nor(0,I_d))$.
\begin{enumerate}[label=(\roman*)]
\item For $q(x)=\kappa_0+b\T x+\tfrac12x\T Ax$ with $A$ symmetric,
\[
\E\big(y-q(x)\big)^2=\big(\E\Phi-\kappa_0-\tfrac12\tr A\big)^2+|g-b|^2+\tfrac12\|\Hb-A\|_F^2+E^{\rm tot}_{\ge3}+\sigma_\xi^2 .
\]
Hence the unrestricted population optimum is $(b,A)=(g,\Hb)$, and under $\operatorname{rank}A\le s$ it is $b=g$, $A=\Hb_s$ (unique if $|\mu_s|>|\mu_{s+1}|$), with value $\mathrm{MSE}_s$.
\item If the current model is $(g,\Hb_s)$, the population PHD matrix of the residual is $\tfrac12(\Hb-\Hb_s)$, whose eigenpair of largest modulus is $(\tfrac12\mu_{s+1},u_{s+1})$. The new column gives $(g,\Hb_{s+1})$, which is optimal within its sign pattern, so ScaledGD does not move it; the validation error falls by $\tfrac12\mu_{s+1}^2=E_2(u_{s+1})$. In population, Step~1 therefore stops at
\[
\hat s=\min\big\{s\ge0:\ \tfrac12\mu_{s+1}^2<\epsilon_0\,\mathrm{MSE}_s\big\}
\]
($\hat s=\operatorname{rank}\Hb$ if the set is empty; in practice also capped at rank 40), and returns the span $\cS_{\rm pop}$ of the eigenvectors of $K_{\hat s}=gg\T+\Hb_{\hat s}^2$ with eigenvalues above $\delta^2$; $\cS_{\rm pop}\subseteq\range(K)\subseteq\cS^*$.
\item For every unit vector $v\perp\cS_{\rm pop}$,
\[
E_1(v)+E_2(v)\le v\T Kv\le\delta^2+\mu_{\hat s+1}^2\le\delta^2+2\epsilon_0\,\mathrm{MSE}_{\hat s},\qquad \mathrm{MSE}_{\hat s}\ge E^{\rm tot}_{\ge3}+\sigma_\xi^2 .
\]
In particular, if $E_2(u_j)<\frac{\epsilon_0}{1-\epsilon_0}\big(E^{\rm tot}_{\ge3}+\sigma_\xi^2\big)$, then $j>\hat s$.
\end{enumerate}
\end{proposition}

The proof is in Appendix~\ref{app:step1}. Once the rank is truncated, Step~1 estimates $\cS_{\rm pop}$ rather than $\range(K)$. By (iii), the resolution is set by the higher-order energy of all directions together, whereas A2 controls it direction by direction; Appendix~\ref{app:step1} gives an example in which A1--A3 hold exactly and Step~1 still misses a well-informed direction. The rule also truncates flat spectra: if the nonzero eigenvalues of $\Hb$ all have the same modulus, the first rank reduces the validation error by at most one over their number, below $\epsilon_0$ when there are more than 50.

For Step~2, the range of $\E M_{3,4}$ is exactly the subspace on which the non-quadratic part of $\Phi$ depends if and only if each of its directions has nonzero third- or fourth-order energy under the design. After $j$ selections the statistic equals the $(j+1)$th largest eigenvalue of $M_{3,4}$, which makes the sequential test asymptotically conservative at the true stopping step; in finite samples it is not exact. If the non-quadratic part within $\range(E)$ is a sum of functions of directions orthogonal in $w$, each with a nonzero third- or fourth-order coefficient, the population maximiser of $\mathcal J$ is unique up to signed permutation, even for equal energies, which the eigenvectors of $M_{3,4}$ cannot achieve. Proofs are in Appendix~\ref{app:step2}.

\subsection{Invariance and mixing}\label{sec:props}

For a $\pi$-reversible kernel $P$ we write $\mathcal E_P(f)=\langle f,(I-P)f\rangle_\pi$ for the Dirichlet form, $\gap(P)=\inf\{\mathcal E_P(f)/\Var_\pi f:\Var_\pi f>0\}$ for the (right) spectral gap, and $\IACT(f)=1+2\sum_{t\ge1}\mathrm{corr}(f(X_0),f(X_t))$ for the integrated autocorrelation time of the stationary chain; the effective sample size per step is $1/\IACT(f)$. Each update leaves $\pi$ (respectively $\pi(x,\tau)$) invariant (Appendix~\ref{app:deriv}): (N) is Metropolis--Hastings on $a$ at fixed $(u,\rho)$ in coordinates with constant Jacobian, (Q) is Metropolis on $c$, the (U) proposals are reversible with respect to the conditional prior, and (H) is Metropolis--Hastings with respect to $\pi(\tau,c,\rho\mid a)$. Errors in the learned structure affect efficiency, not the invariant distribution. Besides the call, an update costs $O(d(r+k)+k^3)$ operations.

For any $\Phi$, in the coordinates $(a,u,\rho)$ the update (N) changes only $a$ (since $c'-\mu(a')=c-\mu(a)$, $u'=u$), (Q) only $u$ ($u'=u+V^{-1/2}S\xi$) and (U) only $\rho$, and each is a Metropolis(--Hastings) kernel reversible with respect to the conditional distribution of its block. The method is thus a random-scan Metropolis-within-Gibbs sampler on $(a,u,\rho)$. All three block kernels are positive operators on $L^2(\pi)$, so the spectrum of the sampler lies in $[0,1]$ and no lazy version is needed (Appendix~\ref{app:exact}).

\label{sec:adapt}After the warm-up, the kernel at step $t$ is $P_t=0.4P_N+0.3P_{Q,S_t}+0.3P_U$ (with a global scale, $0.36P_N+0.27P_{Q,S_t}+0.27P_U+0.1P_H$). Each $P_t$ leaves $\pi$ invariant, but $S_t$ depends on the past, so the chain is adaptive. Ergodicity of adaptive MCMC usually rests on diminishing adaptation and containment \citep{roberts2007,andrieu2008}. Diminishing adaptation always holds here, because the RAM step sizes $\eta_n$ tend to zero (Appendix~\ref{app:Q}). We have not verified containment. \citet{vihola2012} proved a strong law of large numbers for RAM on a fixed target under tail conditions, on the event that the adaptation is stable (eigenvalues of $S_nS_n\T$ bounded and bounded away from zero), which is assumed rather than proved; our (Q) acts on conditional distributions that vary with $(a,\rho)$, and adaptive Metropolis-within-Gibbs samplers generally need further conditions \citep{latuszynski2013}.

\begin{proposition}[exact class]\label{prop:exact}
Suppose $\Phi(x)=\phi_0(a)+(\ell+Ba)\T c+\tfrac12c\T Hc$ exactly, with the sampler's $\ell,B,H$ ($H$ after the eigenvalue limits of Step~3), and $\pi$ normalisable.
\begin{enumerate}[label=(\roman*)]
\item Under $\pi$, $a$, $u$ and $\rho$ are independent, $u\sim\Nor(0,I_k)$, $\rho\sim\Nor(0,P_U)$, and $a$ has density proportional to
\[
m(a)=\exp\big\{-\tfrac12|a|^2-\phi_0(a)+\tfrac12\,\mu(a)\T V^{-1}\mu(a)\big\};
\]
moreover $\alpha_N=\min\{1,m(a')/m(a)\}$ and $\alpha_U\equiv1$.
\item Freeze $S$, $\sigma_j$ and $b_U$. On $L^2(\pi)=L^2(\pi_a)\otimes L^2(\Nor(0,I_k))\otimes L^2(\Nor(0,P_U))$, with $\pi_a\propto m$, the block kernels are $K_N\otimes I\otimes I$, $I\otimes K_Q\otimes I$ and $I\otimes I\otimes K_U$, where $K_N$ is the random-coordinate random-walk Metropolis kernel on $\pi_a$, $K_Q$ the random-walk Metropolis kernel on $\Nor(0,I_k)$ with proposal covariance $\Sigma_u=V^{-1/2}SS\T V^{-\top/2}$, and $K_U$ the equal mixture of refresh and pCN on $\Nor(0,P_U)$, always accepted. With $\gamma_N,\gamma_Q,\gamma_U$ their spectral gaps,
\[
\gap(P)=\min\{0.4\,\gamma_N,\ 0.3\,\gamma_Q,\ 0.3\,\gamma_U\}
\]
($\min\{0.5\gamma_Q,0.5\gamma_U\}$ if $r=0$).
\item $\gamma_U=1-\tfrac12\sqrt{1-b_U^2}$. $\gamma_Q$ depends only on the eigenvalues of $\Sigma_u$ and $\gamma_Q\le4/(ek)$ for every $S$; if $SS\T=\theta V$, then $\Sigma_u=\theta I_k$ and $\gamma_Q$ depends only on $k$ and $\theta$.
\end{enumerate}
\end{proposition}

Thus, in the exact class with $SS\T\propto V$, the spectral gap does not depend on $d$, and the stiffness of the Q directions enters only through the marginal of $a$, a property of the target rather than a cost of the sampler: $\gamma_N$ depends on $\pi_a$, the step sizes and the basis of the N block, $\gamma_Q$ on $k$ and $\theta$, and $\gamma_U\ge\tfrac12$ ($\gamma_U=1$ for $b_U=1$). There is no general formula for $\gamma_N$, which can be exponentially small for multimodal $\pi_a$; by $\gamma_Q\le4/(ek)$, (Q) needs at least about $ek/4$ updates to relax, whatever the step size. $SS\T\propto V$ is the fixed point of RAM in the exact class (Appendix~\ref{app:Q}).

By contrast, for pCN the stiffest directions force a small step and thereby slow down all weakly informed and uninformative directions (Appendix~\ref{app:pcn}).

Learned coordinates are only approximately exact (Section~\ref{sec:learned}). A bounded deviation from the exact class lowers the spectral gap by at most the factor $e^{-\mathrm{osc}\,\zeta}$, where $\mathrm{osc}\,\zeta$ is the oscillation of the deviation $\zeta$ of $\Phi$ (Appendix~\ref{app:exact}); for unbounded deviations, such as errors in the Q coefficients or informative directions left in U, we have no general quantitative result.

In the exact class, ergodic averages over the N and U blocks are unaffected by the adaptation, and those over the Q block reduce to the stability of RAM on a Gaussian target, which to our knowledge is unproved for $k\ge2$ (Appendix~\ref{app:Q}).

\subsection{The global scale}\label{sec:tautheory}

We say that the exact class with a global scale holds if $\Phi(x)=\phi_0(a)+b(a)\T c+\tfrac12c\T H_+c$ for all $x$, with the sampler's own $\ell,B,H_+$; then $R=\phi_0(a)$. Let $\tilde\rho=\rho/\sqrt\tau$ and let $p_N,p_Q,p_U,p_H$ be the selection probabilities.

\begin{proposition}\label{prop:tau}
Assume the exact class with a global scale, a normalisable $\pi$, and a frozen RAM state $S$ and frozen $\sigma_j,b_U$.
\begin{enumerate}[label=(\roman*)]
\item Under $\pi$, $(a,\tau)$, $u$ and $\tilde\rho$ are independent, $u\sim\Nor(0,I_k)$ and $\tilde\rho\sim\Nor(0,P_U)$; $\hat p(\eta\mid a)$ is exactly the conditional posterior density of $\eta$ given $a$, and given $(a,\tau)$, $c\sim\Nor(\mu_\tau(a),V_\tau)$ and $\rho\sim\Nor(0,\tau P_U)$. Hence (H) is a partially collapsed Gibbs update \citep{vandyk2008} with $\alpha_H\equiv1$.
\item $(a_t,\tau_t)$ is a time-homogeneous Markov chain with kernel $M=p_N\mathcal K_N+p_H\mathcal K_H+(1-p_N-p_H)I$, where $\mathcal K_N$ is the (N) kernel on $a$ given $\tau$ and $\mathcal K_H$ the slice kernel on $\tau$ given $a$. $M$ depends neither on $S$, $b_U$, $(u,\tilde\rho)$ nor on $d$.
\end{enumerate}
\end{proposition}

The proof is in Appendix~\ref{app:tauproofs}. Proposition~\ref{prop:tau} makes precise why (H) conditions only on $a$: in the exact class, $\tau$ forms an $(r+1)$-dimensional chain with $a$ while $c$ and $\rho$ are redrawn exactly, so the mixing of $\tau$ does not depend on $d$, and $k$ and the Q stiffness enter only through the marginal of $(a,\tau)$, whatever RAM does. Appendix~\ref{app:tauproofs} also bounds the spectral gap of the whole chain by that of $(a,\tau)$ and gives a dimension-free bound for $r=0$. By contrast, any random-scan sampler that updates $\tau$ only given all of $x$ has a spectral gap of order at most $1/d$ under an inverse-gamma hyperprior with a fixed number of informed directions (Proposition~\ref{prop:centered} in Appendix~\ref{app:tauproofs}).

\section{Numerical experiments}\label{sec:numerics}

\subsection{Synthetic examples}\label{sec:examples}

\begin{figure}[htbp]
\centering
\includegraphics[width=\textwidth]{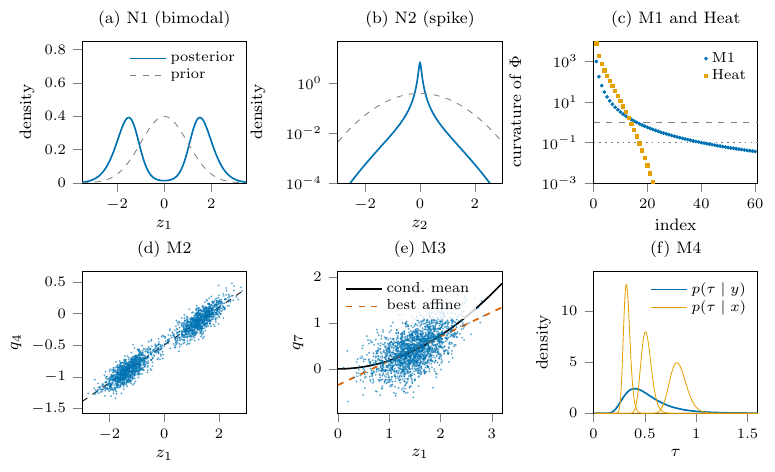}
\caption{Posteriors of the synthetic examples along their characteristic directions. (a) N1 and (b) N2 under the $\Nor(0,1)$ prior (solid; prior dashed). (c) Curvatures of $\Phi$: $\lambda_j$ in M1 and the eigenvalues of the Hessian of $\Phi$ in Heat; dashed: prior precision 1; dotted: $\delta=0.1$. (d) M2: 2000 independent posterior draws of $(z_1,q_4)$ and the conditional mean of $q_4$ (dashed). (e) M3: draws of $(z_1,q_7)$, the conditional mean of $q_7$ and its best affine approximation. (f) M4: posterior density of $\tau$ and the conditional density of $\tau$ given $x$ for three posterior draws of $x$, with $\tau$ at its 10\%, 50\% and 90\% quantiles.}\label{fig:examples}
\end{figure}

The five synthetic examples have $d=200$. Except in M1 and Heat, the likelihood depends on 12 orthonormal directions: two non-Gaussian directions $z_1,z_2$ (N1, N2) and ten Gaussian directions $q_1,\dots,q_{10}$ (Q1--Q10), fixed random projections of $x$; the remaining 188 directions are uninformative. M1--M4 are built from the likelihood
\[
\Phi_G=5e^{-z_1^2}+\log(1+400z_2^2)+\tfrac12\sum_{j=1}^{10}\lambda_j(q_j-m_j)^2,
\]
with $\lambda=(1000,400,200,100,50,25,10,5,3,1.5)$ and $m_j=\pm0.5$ alternating. Under the $\Nor(0,1)$ prior, N1 is bimodal with modes at $\pm1.52$ and a barrier of 3.35 nats, and N2 is a spike with posterior standard deviation 0.20 and kurtosis 27. Table~\ref{tab:toys} defines the examples and Fig.~\ref{fig:examples} shows their posteriors along the characteristic directions.

\begin{table}[htbp]
\centering\small
\caption{Synthetic examples ($d=200$). M1--M3 have prior $\Nor(0,I_{200})$.}\label{tab:toys}
\begin{tabularx}{\textwidth}{>{\raggedright\arraybackslash}p{2.1cm}>{\raggedright\arraybackslash}X>{\raggedright\arraybackslash}p{3.3cm}}
\toprule
Example & Definition & What it tests \\
\midrule
M1 (many weakly informed directions) & N1, N2 as in $\Phi_G$; 190 Gaussian directions with $\lambda_j=1000\,j^{-2.5}$ (175 weakly informed) and $m_j=\pm0.5$; 8 uninformative directions & separating informed from uninformative directions when there is no clear cutoff \\ \addlinespace
M2 (coupling with a bimodal direction) & $\Phi_G$ with $q_j-m_j$ replaced by $q_j-m_j-0.3z_1$, so that the conditional means of the $q_j$ move linearly with N1 & jumping between modes when the Gaussian directions must move with the bimodal one \\ \addlinespace
M3 (curved posterior) & $\Phi_G$ with the N1 term replaced by $\tfrac12\cdot3.605\,(z_1-1.938)^2$ (unimodal) and $q_j-m_j$ by $q_j-m_j-b_j(z_1^2-c_0)$, with $b_j=0.2$ for Q7--Q10, $b_j=0$ otherwise and $c_0=2.52$; banana-shaped joint posterior of N1 and Q7--Q10 & robustness to a curved posterior \\ \addlinespace
M4 (non-Gaussian prior with a global scale) & $\theta\mid\tau\sim\Nor(0,\tau I_{200})$, $\tau\sim\mathrm{InvGamma}(3,2)$ ($\E\tau=1$, sd 1); likelihood $\Phi_G$ & handling a non-Gaussian prior together with non-Gaussian directions \\ \addlinespace
Heat (backward heat conduction problem) & Initial condition $u_0(s)=\sum_{k=1}^{200}(2\lambda_k)^{1/2}x_k\sin(k\pi s)$, $\lambda_k=(k\pi)^{-2}$, of $u_t=u_{ss}$ on $(0,1)$ with zero boundary values, observed as $y_j=u(j/31,T)+\epsilon_j$, $j=1,\dots,30$, $T=10^{-3}$, $\epsilon_j\sim\Nor(0,0.02^2)$; prior as in M4, data generated with $\tau=0.5$; 16 eigenvalues of the Hessian of $\Phi$ exceed 0.1 & a control case in which all directions are Gaussian \\
\bottomrule
\end{tabularx}
\end{table}

Reference values of the posterior means and covariances of all coordinates and of $\tau$ come from one- or two-dimensional quadrature, over the non-Gaussian directions and, where present, $\log\tau$, combined with the Gaussian conditionals of the other directions. The posterior of $\tau$ has mean 0.528 and 90\% interval 0.26--0.95 in M4, and mean 0.781 and 90\% interval 0.41--1.35 in Heat. All N and Q directions satisfy A2 except the spike N2, which is on the boundary, and A3 holds exactly except in M3 (Appendix~\ref{app:energy}).

We compare the NQU sampler with nine baselines (Appendix~\ref{app:exp}). The gradient-free ones are full-space pCN \citep{cotter2013}; gpCN \citep{rudolf2018}, with the empirical mean and covariance of the chain; QPHD, which estimates the informed subspace from a full quadratic model \citep{constantine2017}; and SAFES and SAFES-P \citep{dunlop2022}. The gradient-based ones are DILI \citep{cui2016}, which also uses Hessians; NUTS \citep{hoffman2014} with a diagonal or a dense mass matrix; and AS-MCMC \citep{constantine2016}, made exact by sampling the inactive variables by MCMC. AS-MCMC and QPHD use the sampling stage of the NQU sampler without an N block, so the comparison with them isolates the effect of the N block. Each chain has a budget of $5\times10^5$ calls, which includes learning, adaptation and warm-up; a gradient, a Hessian and a leapfrog step of NUTS each count as one call, which favours the gradient-based methods. Dense NUTS exceeded the budget in its warm-up alone and is reported without it. With a global scale, pCN and NUTS use a non-centred parameterisation, AS-MCMC and QPHD the (H) update, and DILI, gpCN, SAFES and SAFES-P, which assume a Gaussian prior, a non-centred parameter with a standard normal prior (Appendix~\ref{app:exp}). The NQU sampler uses $n_1=30\,000$, $n_2=20\,000$ and $4.5\times10^5$ sampling calls, the first 10\% as warm-up; we also report $n_1=20\,000$, $n_2=30\,000$.

Following \citet{hoffman2014}, efficiency is measured by the minimum, over the means and second central moments of all coordinates, of the effective sample size (ESS) per 1000 calls, estimated by Geyer's initial monotone sequence \citep{geyer1992}. It is computed in two coordinate systems: on the posterior axes, whose leading axes are the informed directions of the example, so that each difficult direction is measured on its own, and in the original coordinates, each of which mixes all directions; with a global scale both include $\tau$. We report the calls needed for 1000 effective samples, that is, the calls before the ESS window (learning, adaptation and warm-up) plus $10^6$ divided by the minimum ESS per 1000 calls, as medians over 8 chains (seeds 100--107). A chain is flagged as biased if $\max|z|>4$ over all quantities of a coordinate system, with $z=(\bar f-f_{\rm true})/(\hat s_f/\mathrm{ESS}_f^{1/2})$; for independent standard normal $z$ this happens with probability about 2.5\% per chain.

\label{sec:results}\label{sec:learned}Tables~\ref{tab:learn} and~\ref{tab:main} summarise the results. Learning recovered the intended structure. All seeds of M1--M4 select $r=2$ and put at least 0.95 of N2 in the N block. In M1 (many weakly informed directions), the 20--22 most strongly informed of the 190 Gaussian directions enter $\cS$, and about 3 nats of information remain in U. In M2 (coupling with a bimodal direction), only 0.56--0.58 of N1 lies in N, as expected, because the back-map of Step~2 includes in the N direction the movement of the coupled Q directions. In M3 (curved posterior), N1 is unimodal and its non-Gaussianity comes only from the coupling, so the selected N1 direction varies between seeds. In Heat, no non-Gaussian direction is selected, and $\cS$ captures at least 0.999 of each of the 16 leading informed directions. Appendix~\ref{app:exp} gives further diagnostics of the learning.

\begin{table}[htbp]
\centering\small
\caption{Learned structures ($n_1=30\,000$, $n_2=20\,000$; ranges over 8 seeds). ``In N'' is the squared length of the projection of the true unit direction on $\range(W_N)$, similarly for Q; the information left in U is $\tfrac12\sum_j\log(1+\lambda_j|P_Uv_j|^2)$ over the Gaussian directions $v_j$.}\label{tab:learn}
\begin{tabularx}{\textwidth}{>{\raggedright\arraybackslash}p{0.9cm}>{\raggedright\arraybackslash}p{1.1cm}>{\raggedright\arraybackslash}p{0.5cm}>{\raggedright\arraybackslash}X>{\raggedright\arraybackslash}p{1.6cm}>{\raggedright\arraybackslash}X}
\toprule
& $\dim\cS$ & $r$ & N1 & N2 & Gaussian directions; information left in U \\
\midrule
M1 & 23--25 & 2 & in N 0.94--0.98 & in N 0.96--0.98 & 20--22 of 190 more than half in $\cS$; 2.97--3.29 nats \\
M2 & 13 & 2 & in N 0.56--0.58, in Q 0.41--0.42 & in N 0.99 & all 10 more than half in $\cS$, at most 0.05 in N; 0.84--0.99 nats \\
M3 & 13 & 2 & in N 0.33--0.91, in Q 0.08--0.66 & in N 0.95--0.96 & all 10 more than half in $\cS$, at most 0.065 in N; 0.10--0.11 nats \\
M4 & 13 & 2 & in N 0.97--0.99 & in N 0.98--0.99 & all 10 more than half in $\cS$, at most 0.005 in N; 0.01--0.02 nats \\
Heat & 16 & 0 & --- & --- & 16 leading directions at least 0.999 in $\cS$; 0.078 nats \\
\bottomrule
\end{tabularx}
\end{table}

\begin{table}[htbp]
\centering\footnotesize\setlength{\tabcolsep}{3pt}
\caption{Calls needed for 1000 effective samples relative to the NQU sampler ($n_1=30\,000$, $n_2=20\,000$), as ratios of medians over 8 chains, posterior axes / original coordinates; values above 1 favour the NQU sampler. The first row gives the NQU calls in units of $10^5$. For each example and coordinate system, the fewest calls are in bold and the second fewest underlined; the row NQU (20k/30k) is not ranked.}\label{tab:main}
\begin{tabularx}{\textwidth}{>{\raggedright\arraybackslash}p{3.1cm}*{5}{>{\raggedright\arraybackslash}X}}
\toprule
Method & M1 & M2 & M3 & M4 & Heat \\
 & \emph{weakly informed} & \emph{coupling with a bimodal direction} & \emph{curved posterior} & \emph{non-Gaussian prior with a global scale} & \emph{backward heat conduction} \\
\midrule
NQU (30k/20k), $\times10^5$ & \textbf{4.44} / \textbf{1.54} & \textbf{2.25} / \textbf{1.19} & \underline{2.38} / 1.34 & \textbf{2.14} / \underline{1.57} & \textbf{1.21} / \textbf{1.21} \\
NQU (20k/30k) & 0.98 / 1.00 & 0.96 / 1.01 & 1.01 / 1.07 & 0.96 / 1.05 & 1.00 / 1.00 \\
\midrule
\multicolumn{6}{l}{\emph{Gradient-based}} \\
DILI & \underline{1.79} / \underline{1.24} & \underline{1.60} / \underline{1.22} & \textbf{0.72} / \underline{0.85} & \underline{1.25} / \textbf{0.86} & 2.20 / 2.02 \\
NUTS (diagonal) & 17 / 5.53 & 14 / 5.39 & 3.24 / 4.46 & 9.93 / 3.91 & 2.44 / 2.44 \\
NUTS (dense) & 4.36 / 6.27 & 5.89 / 7.64 & 4.03 / 6.72 & 5.47 / 5.89 & 6.59 / 6.59 \\
AS-MCMC (exact) & 16 / 5.20 & 6.34 / 2.08 & 5.32 / \textbf{0.75} & 4.69 / 6.40 & $\infty$\footnotemark[1] \\
\midrule
\multicolumn{6}{l}{\emph{Gradient-free}} \\
pCN & 57 / 38 & 157 / 96 & 36 / 71 & 82 / 39 & 355 / 349 \\
gpCN & 145 / 404 & 272 / 504 & 256 / 468 & 283 / 377 & 511 / 526 \\
QPHD & 48 / 90 & 78 / 68 & 155 / 203 & 11 / 15 & \underline{1.01} / \underline{1.01} \\
SAFES & 18 / 23 & 34 / 28 & 35 / 26 & 33 / 20 & 51 / 51 \\
SAFES-P & 18 / 23 & 37 / 26 & 37 / 24 & 22 / 15 & 23 / 23 \\
\bottomrule
\end{tabularx}
\footnotetext[1]{The (H) update of AS-MCMC was never accepted, so $\tau$ kept its initial value in all 8 chains (Appendix~\ref{app:exp}).}
\end{table}

\begin{figure}[htbp]
\centering
\includegraphics[width=\textwidth]{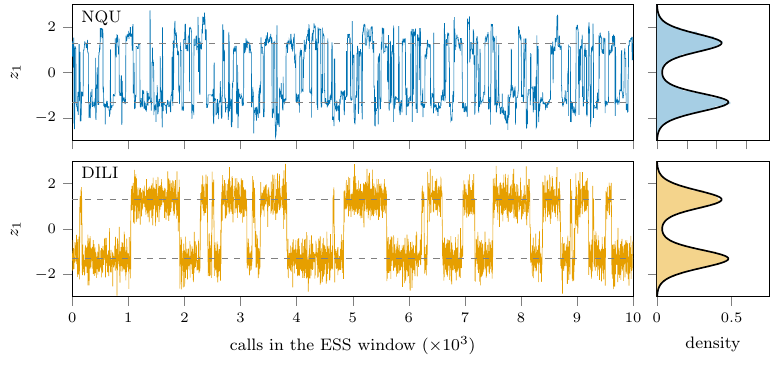}
\caption{M2: coordinate $z_1$ along the true N1 direction in one chain of the NQU sampler and one of DILI (seed 100), over the first $10^4$ calls of the ESS window; dashed lines mark the two modes. Right: histograms over the whole window ($4\times10^5$ calls) with the exact marginal (solid). Over the whole window, $z_1$ moves between the modes 6773 times in the NQU chain and 1481 times in the DILI chain.}\label{fig:traces}
\end{figure}

On the posterior axes, the NQU sampler needed the fewest calls in M1, M2, M4 and Heat; in the original coordinates it did so in M1, M2 and Heat. In M1 and M2 the next best method is DILI, with 1.79 and 1.60 times as many calls on the posterior axes and 1.24 and 1.22 times in the original coordinates; the gap is smaller in the original coordinates because both methods spend about $9.5\times10^4$ calls before the ESS window, and within the window the NQU sampler is 1.6 and 2.1 times faster. In all 16 chains the bottleneck of DILI is the mean of the bimodal N1: its Langevin-type proposal crosses the barrier between the modes less often than the random walk of (N), which is accepted at rate 0.34--0.36 (Fig.~\ref{fig:traces}). AS-MCMC, which samples like the NQU sampler but without an N block, needed 4.7--16 times the NQU calls on the posterior axes in M1, M2 and M4, probably because one Gaussian proposal then covers N1, N2 and the Gaussian directions. In M4 (non-Gaussian prior with a global scale) the next best method on the posterior axes is DILI, applied to the non-centred parameter, with 1.25 times the NQU calls. In the original coordinates, which include $\tau$, DILI needed only 0.86 times the NQU calls: there the worst quantity of the NQU sampler is $\tau$ in all chains, and DILI samples $\tau$ with a median effective sample size per call 1.7 times as large. AS-MCMC and QPHD use the same (H) update as the NQU sampler but accept it at rate 0.02--0.10, against 0.56--0.59. In Heat, QPHD needed 1.01 times the NQU calls in both coordinate systems: for this linear problem its quadratic model is exact, and with $r=0$ its sampling stage is that of the NQU sampler. DILI needed 2.20 and 2.02 times the NQU calls, the two NUTS variants 2.44 and 6.59 times, and AS-MCMC never updated $\tau$, because its Gaussian model from Stein's identity is inaccurate in the stiffest directions (Appendix~\ref{app:exp}). The bottleneck of the NQU sampler there is the scale $\tau$, and property (iv) in Appendix~\ref{app:tauproofs} predicts its efficiency: 47.2 and 39.0 effective samples per 1000 calls for $\tau$ and $(\tau-\E\tau)^2$ in the exact class, against observed medians of 45.75 and 38.95 (the learned structure has $\alpha_H=0.95$ instead of 1); DILI attains 18.7 and 7.1.

M3 (curved posterior) is the one example in which the NQU sampler is the fastest in neither coordinate system: DILI needed 0.72 and 0.85 times its calls, and AS-MCMC 0.75 times in the original coordinates, although 5.32 times on the posterior axes. In 6 of 8 chains the worst quantity of the NQU sampler on the posterior axes is the mean of Q7, the direction with the most curved conditional mean, which is consistent with the violation of A3: the Q model is affine in $a$ and cannot follow the quadratic dependence of Q7 on N1 (Appendix~\ref{app:exp}).

Apart from QPHD in Heat, the gradient-free baselines needed at least 10.8 times the NQU calls, and at least 17.7 times in M1--M3, where most of their chains are biased; the two NUTS variants needed 2.4--17 times the NQU calls (Appendix~\ref{app:exp}). Of the 80 chain--coordinate-system pairs of the NQU sampler, three are flagged as biased, with $\max|z|$ between 4.13 and 4.59, slightly more than the nominal 2.5\%. The two learning budgets give similar results (ratios 0.96--1.07), and learning accounts for 11--42\% of the calls needed for 1000 effective samples.

\subsection{Application: reaction rates in methane ignition}\label{sec:cantera}

This example uses a real reaction mechanism with simulated data, so no reference values are available. We simulate the ignition of a stoichiometric CH$_4$/air mixture in a constant-volume zero-dimensional reactor with Cantera 3.2 \citep{cantera} and GRI-Mech 3.0 \citep{grimech} without nitrogen chemistry (36 species, 219 reactions). The unknowns are the log-multipliers of all 219 rate constants, $\ln k_i=\ln k_i^0+\sigma_ix_i$ with $x\sim\Nor(0,I_{219})$, so $d=219$; two reactions, C$_2$H$_6$+O $\rightleftharpoons$ C$_2$H$_5$+OH and CH$_3$+O $\to$ CO+H+H$_2$, have a wide prior, $\sigma_i=2.0$, and all others have $\sigma_i=0.25$. At 9 initial conditions (temperature 1150, 1400 and 1800\,K, pressure 1, 10 and 30\,atm) we observe the logarithms of the ignition delay and of the peak mole fractions of CH$_2$O and C$_2$H$_6$, 27 observations in total, simulated at $x=0$ with Gaussian noise of standard deviation 0.20. One call integrates 9 reactors and takes about 0.4\,s. Cantera's zero-dimensional reactor has no adjoint, and forward sensitivities with respect to 200 rate multipliers cost as much as 192 plain solves, so we treat the model as a black box; further details are in Appendix~\ref{app:exp}.

\begin{figure}[htbp]
\centering
\includegraphics[width=\textwidth]{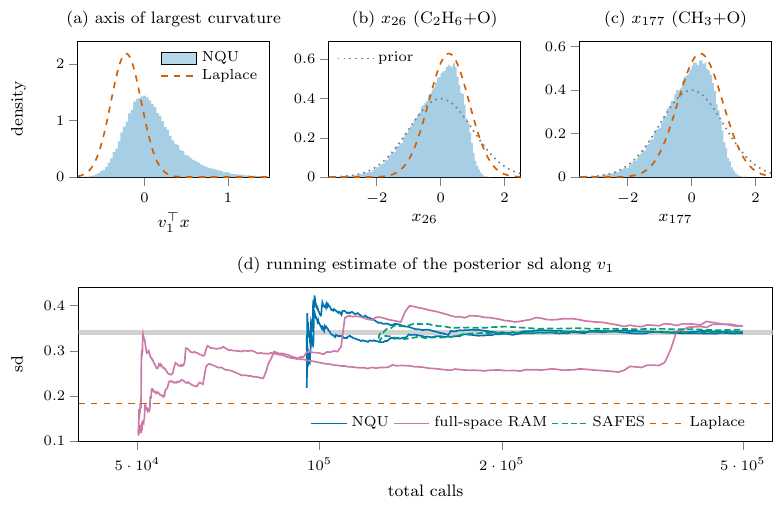}
\caption{Methane example. (a)--(c) Histograms of the two pooled NQU chains along the Laplace axis $v_1$ of largest curvature and for the two reactions with a wide prior, with the Laplace approximation (dashed) and, in (b) and (c), the prior (dotted). (d) Running estimates of the posterior standard deviation along $v_1$ against total calls, two chains per method, each from the start of its ESS window; grey band: pooled NQU estimate $\pm2$ standard errors; dashed: Laplace value. Along $v_1$ the SAFES estimates are close to the NQU one; its disagreement concerns other quantities (Table~\ref{tab:cantera}).}\label{fig:cantera}
\end{figure}

The Laplace approximation has 5 curvatures above 1 at the mode, the other directions being close to the prior. According to the NQU chains, the posterior is clearly non-Gaussian: along the direction of the largest curvature it has skewness about $+1.05$ and a standard deviation about 1.9 times that of the Laplace approximation (Fig.~\ref{fig:cantera}a). The two reactions with a wide prior are constrained from one side only, since a large rate contradicts the data while a small one hardly affects them, and the posterior along them has skewness about $-0.7$ (Fig.~\ref{fig:cantera}b, c).

We run the NQU sampler as in Section~\ref{sec:examples}, with one chain for each of the seeds 100 and 101, and compare it with two gradient-free baselines with the same seeds and $5\times10^5$ calls per chain: SAFES, and full-space RAM \citep{vihola2012}, the adaptive kernel of the (Q) update applied to all 219 coordinates, which amounts to the NQU sampler without learning. ESS is computed as in Section~\ref{sec:examples} for 438 quantities per coordinate system, with the eigenvectors of the Laplace covariance as posterior axes. Without reference values we report agreement: the number of quantities with $|z|>4$ between the two seeds of a method, and between each chain and the pooled estimate of the two NQU chains.

\begin{table}[htbp]
\centering\footnotesize\setlength{\tabcolsep}{3pt}
\caption{Methane example, $5\times10^5$ calls per chain. Each cell gives posterior axes / original coordinates; pairs of numbers are seeds 100 and 101. ``Against NQU'' counts quantities (of 438) with $|z|>4$ between a chain of the method and the pooled estimate of the two NQU chains. In the first three rows, the best values for each coordinate system are in bold.}\label{tab:cantera}
\begin{tabular}{>{\raggedright\arraybackslash}p{2.5cm}lll}
\toprule
 & NQU & SAFES\footnotemark[1] & Full-space RAM \\
\midrule
Minimum ESS per 1000 calls & \textbf{14.1, 13.6} / \textbf{12.5, 12.8} & 0.48, 0.46 / 0.43, 0.36 & 0.020, 0.020 / 0.018, 0.025 \\
Calls for 1000 ESS ($\times10^5$) & \textbf{1.66, 1.69} / \textbf{1.75, 1.73} & 21.9, 23.0 / 24.3, 28.9 & 512, 510 / 551, 398 \\
$|z|>4$ between seeds & \textbf{0} / \textbf{0} & 57 / 66 & 3 / 4 \\
$|z|>4$ against NQU & --- & 72, 74 / 74, 70 & 11, 5 / 12, 8 \\
\bottomrule
\end{tabular}
\footnotetext[1]{SAFES degenerated: after its step parameter fell to zero, its particles could not leave the affine hull of the ensemble, so its ESS measures mixing within that hull only and does not reflect the accuracy of its estimates (Appendix~\ref{app:exp}).}
\end{table}

Table~\ref{tab:cantera} gives the results. Both NQU seeds learned the same structure: no N direction ($r=0$), a four-dimensional Q block and 215 directions in U. Their estimates agree ($\hat R\le1.0006$ and $|z|\le3.5$ for all 438 quantities), and they need $1.66$--$1.75\times10^5$ calls for 1000 effective samples. Step~2 selected no N direction because its design is centred at the Step~1 estimate, where $-\log p$ is about 21 nats higher than near the posterior and the conditionals along the two reactions with a wide prior are nearly Gaussian; the skewed directions are therefore sampled by the Q and U updates. Full-space RAM largely agrees with the NQU sampler (3 and 4 of the 438 quantities with $|z|>4$ for its pooled chains) but mixes poorly within the budget (Fig.~\ref{fig:cantera}d), with $\hat R$ up to 1.19, and needs about 300 times the NQU calls; the same RAM kernel acts on the learned four-dimensional Q block in one case and on all 219 coordinates in the other. SAFES degenerated: its tuning rule drove its step parameter to zero, which confined the particles to the affine hull of the ensemble and froze 180 directions, and its two seeds disagree in 57 and 66 of the 438 quantities (Appendix~\ref{app:exp}). Of the three methods, only the NQU sampler produced estimates on which its two seeds agree within $5\times10^5$ calls; full-space RAM largely agrees with it, at a cost two orders of magnitude higher.

\section{Conclusion}\label{sec:conclusion}

We have proposed the NQU sampler, a gradient-free MCMC method for Bayesian inverse problems whose posteriors are approximately Gaussian in most informed directions but strongly non-Gaussian in a few. It learns a split into non-Gaussian, approximately quadratic and uninformative directions from likelihood values and uses the split only to construct proposals, so each update leaves the exact posterior invariant at the cost of one call. For likelihoods with exactly the learned structure, its spectral gap does not depend on the dimension, and neither does the mixing of a global prior scale under the collapsed update. In the experiments, with learning included, it needed fewer calls than gradient-based baselines in most cases and, in the examples with non-Gaussian directions, at least an order of magnitude fewer than gradient-free ones.

The method has several limitations. Learning relies on A2, which cannot be checked before learning, and even when A2 holds, the higher-order energy of all directions together can make Step~1 miss an informed direction (Section~\ref{sec:step1pop}). A missed direction is left in U, which does not change the invariant distribution, but a strongly informed one, such as a spike, can make most refreshes of U fail; this did not happen in the examples of Section~\ref{sec:numerics}. The design of Step~2 is centred at the Step~1 estimate, so non-Gaussian directions can go undetected when the posterior is concentrated away from it, as in the methane example, where no N direction was selected. The conditional mean of the Q model is affine in the N coordinates, so curved coupling as in M3 is followed only approximately. The mixing results cover only the exact class and bounded deviations from it, and the containment condition of the adaptation is not verified. Finally, the synthetic examples are exact ridge functions, which favours all subspace methods, $d$ is fixed at 200 there, and the application has no reference values and two chains per method.

These limitations suggest further work. Adaptive or sequential designs for Step~2 that move towards the posterior could find non-Gaussian directions near a posterior concentrated away from the Step~1 estimate, and a Q model with conditional means nonlinear in the N coordinates could follow curved coupling. How the learning budget must grow with $d$ also remains to be studied: in finite samples the stopping threshold of Step~1 increases with $d/n_1$ (Appendix~\ref{app:exp}), so a fixed $n_1$ may not suffice for larger $d$.

\backmatter

\section*{Statements and Declarations}

\begin{itemize}
\item Competing interests: the authors declare no competing interests.
\item Funding: Jiang and Liu were partially supported by the Henan Provincial International Science and Technology Cooperation Project (No.~252102521018).
\item Data availability: all data used in this paper are synthetic.
\item Code availability: a Python implementation of the NQU sampler is available at \url{https://github.com/lijinglai/nqu-sampler}.
\end{itemize}

\begin{appendices}

\section{Proofs}\label{app:deriv}

\subsection{Back-map and exactness}\label{app:backmap}
Whether the learned coordinates are exact depends on the back-map of Step~2. Let $\Phi=q+\psi(\Gamma\T x)$, where $q$ is the quadratic part of $\Phi$ and $\psi$ depends on $r$ directions $\Gamma=[\gamma_1,\dots,\gamma_r]$ in $\cS$ and is not quadratic along any of them. For any orthonormal $W_N$ in $\cS$ with complement $W_Q$, if the floor on $I_k+H$ is inactive, completing the square as in Appendix~\ref{app:exact} gives
\[
\pi(a,u,\rho)\ \propto\ m_Q(a)\,e^{-|u|^2/2}e^{-|\rho|^2/2}\exp\{-\psi(C_0+C_1a+N_uu)\},\qquad N_u=\Gamma\T W_QV^{1/2},
\]
with $m_Q$ determined by $q$ and constant $C_0,C_1$. The learned coordinates are therefore exact if and only if $N_u=0$, that is $\range(\Gamma)\subseteq\range(W_N)$. When Step~2 is identifiable, $E$ spans $A_{\cS}^{-1/2}F_{\cS}\T\range(\Gamma)$ in population, so $W_N$ spans $F_{\cS}A_{\cS}^{-1}F_{\cS}\T\range(\Gamma)$, which contains $\range(\Gamma)$ if and only if $F_{\cS}\T\range(\Gamma)$ is invariant under $A_{\cS}$, equivalently $B=0$ when the floor on $A_{\cS}$ is inactive; in M2, where the floor is active, $N_u\ne0$ as well. With $W_N=\orth(F_{\cS}A_{\cS}^{1/2}E)$ one would have $N_u=0$ in population. The deviation is measured by $N_u$: $\pi(u\mid a)$ differs from $\Nor(0,I_k)$ only through $N_uu$, and for $r=1$, $|N_u|$ is the standard deviation of $\gamma_1\T x$ given $a$ under the Q model; it is about 0.11 for the population structure of M2, where the deviation is mainly a shift of conditional mean and variance and $\pi(u\mid a)$ remains close to Gaussian.

\subsection{Acceptance probability of (H)}\label{app:H}
Fix $a$. Let $q_\tau(c)=\Nor(c;\mu_\tau(a),V_\tau)$ and $g_\tau$ the density of $\Nor(0,\tau P_U)$ on $\range(P_U)$. Given $\tau$ the blocks are independent a priori, and completing the square gives $\Nor(c;0,\tau I_k)e^{-\Phi_Q(a,c)}=Z_Q(a;\tau)\,q_\tau(c)$, so
\[
\pi(\tau,c,\rho\mid a)\ \propto\ \hat p(\tau\mid a)\,q_\tau(c)\,g_\tau(\rho)\,e^{-R(x)},\qquad \hat p(\tau\mid a)\propto p(\tau)\Nor(a;0,\tau I_r)Z_Q(a;\tau)
\]
(with Jacobian $e^\eta$ in $\eta=\log\tau$). Let $\mu_a(d\tau,dc,d\rho)=\hat p(\tau\mid a)q_\tau(c)g_\tau(\rho)\,d\tau\,dc\,d\rho$, the possibly only $\sigma$-finite distribution obtained by replacing $\Phi$ with $\Phi_Q$; $\pi(\cdot\mid a)$ has density $h\propto e^{-R}$ with respect to $\mu_a$, and $\int h\,d\mu_a<\infty$. The proposal is $Q_a(y,dy')=K_a(\tau,d\tau')\,q_{\tau'}(c')\,g_{\tau'}(\rho')\,dc'\,d\rho'$ with $y=(\tau,c,\rho)$ and $K_a$ the slice kernel, which is reversible with respect to $\hat p(\tau\mid a)\,d\tau$ even if $\hat p$ is not integrable, since the argument of \citet{neal2003} uses only the symmetry of Lebesgue measure on the slice. Hence $\mu_a(dy)Q_a(y,dy')$ is symmetric in $(y,y')$, and so is $\pi(dy\mid a)Q_a(y,dy')\min\{1,h(y')/h(y)\}\propto\mu_a(dy)Q_a(y,dy')\min\{h(y),h(y')\}$, which is detailed balance; $h(y')/h(y)=e^{R(x)-R(x')}$. (H) does not change $a$, so it also leaves $\pi(x,\tau)$ invariant.

\subsection{Integrability of $\hat p$}\label{app:phat}
The density $\hat p(\cdot\mid a)$ used in (H) need not be integrable. Write $H_+=\sum_ih_ie_ie_i\T$ and let $P_0$ be the orthogonal projection onto the null space of $H_+$. Then
\[
Z_Q(a;\tau)=\prod_{h_i>0}(1+\tau h_i)^{-1/2}\exp\Big\{\sum_{h_i>0}\frac{\tau\,(e_i\T b(a))^2}{2(1+\tau h_i)}+\frac\tau2\,|P_0b(a)|^2\Big\}.
\]
If $H_+$ has a zero eigenvalue, for instance because a negative eigenvalue of $H$ was set to zero, and $P_0b(a)\ne0$, then $Z_Q(a;\tau)$ grows exponentially in $\tau$, and for a hyperprior with polynomial tails, such as the inverse gamma, $\hat p(\cdot\mid a)$ is not integrable. This occurs in M4 (Section~\ref{sec:numerics}); for seed 100, an eigenvalue $-0.0034$ of $H$ is set to zero. Correctness is unaffected, because Appendix~\ref{app:H} uses only the reversibility of the slice kernel with respect to the $\sigma$-finite measure $\hat p(\tau\mid a)\,d\tau$. Efficiency may suffer: the slice in $\eta$ is then unbounded above, and proposals at large $\tau$ are rejected in step (3) of (H).

\subsection{Proof of Proposition~\ref{prop:K}}\label{app:K}
The first two chaos components are $\Phi_1(x)=g\T x$ and $\Phi_2(x)=\tfrac12\sum_{ij}\Hb_{ij}(x_ix_j-\delta_{ij})$, with $\|\Phi_2\|^2=\tfrac12\|\Hb\|_F^2$. In an orthonormal basis whose first vector is $v$, let $\tilde H$ represent $\Hb$. The part of $\Phi_1$ involving $v$ is $(v\T g)(v\T x)$, so $E_1(v)=(v\T g)^2$, and the part of $\Phi_2$ not involving $v$ is the quadratic form of $\tilde H$ without its first row and column, so
\[
E_2(v)=\tfrac12\Big(\|\tilde H\|_F^2-\sum_{i,j\ge2}\tilde H_{ij}^2\Big)=\tfrac12\Big(\tilde H_{11}^2+2\sum_{j\ge2}\tilde H_{1j}^2\Big)=|\Hb v|^2-\tfrac12(v\T\Hb v)^2 ,
\]
which is \eqref{eq:E12}. (i) $v\T Kv=(v\T g)^2+|\Hb v|^2$ and $0\le(v\T\Hb v)^2\le|\Hb v|^2$. (ii) Under A2, $E_{\rm tot}\le2(E_1+E_2)\le2v\T Kv$; and $\E[\Phi\mid P_{v^\perp}x]$ is the sum of the terms not involving $v$, so $E_{\rm tot}(v)=\E[\Var(\Phi\mid P_{v^\perp}x)]$. (iii) Let $P=P_{\cS^*}$. Since $\Phi$ depends only on $Px$, which is independent of $(I-P)x$, $(I-P)g=0$ and $\Hb u=0$ for $u\perp\cS^*$; so $g\in\cS^*$, $\Hb$ maps $\cS^*$ into itself and vanishes on its complement, and $\range(K)\subseteq\cS^*$. As $K\succeq0$, $\cS^*\subseteq\range(K)$ if and only if $\cS^*\cap\ker K=\{0\}$, and by (i) $Kv=0$ if and only if $E_1(v)=E_2(v)=0$; if $\Hb$ is nonsingular on $\cS^*$, then $v\T Kv\ge|\Hb v|^2>0$ for $0\ne v\in\cS^*$.

\subsection{Proof of Proposition~\ref{prop:step1}}\label{app:step1}
(i) Write $q=(\kappa_0+\tfrac12\tr A)+b\T x+\tfrac12\langle A,xx\T-I\rangle$. The chaos components of orders 0--2 of $\Phi$ are $\E\Phi$, $g\T x$ and $\tfrac12\langle\Hb,xx\T-I\rangle$, and by Isserlis' formula $\E[\langle A,xx\T-I\rangle\langle C,xx\T-I\rangle]=2\langle A,C\rangle_F$ for symmetric $A,C$; orthogonality of the chaos and $\E[\xi\mid x]=0$ give the decomposition. Singular values of a symmetric matrix are the moduli of its eigenvalues and its truncated singular value decomposition is $\Hb_s$, so by the Eckart--Young--Mirsky theorem $\Hb_s$ is the best Frobenius approximation of rank at most $s$; $LDL\T$ with free signs ranges over all symmetric matrices of rank at most $s$.

(ii) For the model $(g,\Hb_s)$ with optimal intercept, the residual is $\tfrac12\langle\Hb-\Hb_s,xx\T-I\rangle+\sum_{k\ge3}\Phi_k+\xi$; chaos of order three or more and $\xi$ are orthogonal to the entries of $xx\T-I$, so the population PHD matrix is $\tfrac12(\Hb-\Hb_s)$. The new column $\sqrt{|\mu_{s+1}|}\,u_{s+1}$ with sign $\operatorname{sign}\mu_{s+1}$ makes $LDL\T=\Hb_{s+1}$, the optimum of rank at most $s+1$ by (i) and hence within the current sign pattern, so the gradient vanishes. The validation error falls by $\mathrm{MSE}_s-\mathrm{MSE}_{s+1}=\tfrac12\mu_{s+1}^2=E_2(u_{s+1})$ by \eqref{eq:E12}, and induction from rank 0 gives $\hat s$. Finally $\range(K_{\hat s})=\mathrm{span}(g)+\range(\Hb_{\hat s})\subseteq\range(K)$, and Proposition~\ref{prop:K}(iii) applies.

(iii) If $v\perp\cS_{\rm pop}$, then $v$ lies in the eigenspace of $K_{\hat s}$ for eigenvalues at most $\delta^2$, so $v\T K_{\hat s}v\le\delta^2$. Since $\Hb_{\hat s}$ and $\Hb-\Hb_{\hat s}$ have orthogonal ranges, $\Hb^2=\Hb_{\hat s}^2+(\Hb-\Hb_{\hat s})^2$ and $v\T Kv\le\delta^2+\mu_{\hat s+1}^2$. The left inequality is Proposition~\ref{prop:K}(i), and the right one follows from the stopping condition (with $\mu_{\hat s+1}=0$ if $\hat s=\operatorname{rank}\Hb$). If the greedy path reaches step $j$, then $\mathrm{MSE}_{j-1}\ge\tfrac12\mu_j^2+E^{\rm tot}_{\ge3}+\sigma_\xi^2$, and the stated condition is equivalent to $\tfrac12\mu_j^2<\epsilon_0(\tfrac12\mu_j^2+E^{\rm tot}_{\ge3}+\sigma_\xi^2)\le\epsilon_0\mathrm{MSE}_{j-1}$, so the rule stops at $j-1$; otherwise $\hat s<j-1$.

The following example shows that Step~1 can miss a direction even when A1--A3 hold exactly. Let $q_1,q_2,z$ be the coordinates of $x$ along three orthonormal directions, $0<\omega\le\lambda/(2\sqrt2)$, $0<\nu<\lambda$, and
\[
\Phi=\tfrac\lambda2\,\mathrm{He}_2(q_1)+\tfrac\lambda2\,\mathrm{He}_2(q_2)+\omega\,\mathrm{He}_2(q_1)\,\mathrm{He}_2(q_2)+\tfrac\nu2\,\mathrm{He}_2(z).
\]
Then $\Phi$ is bounded below and, with N block $(q_1,q_2)$ and Q block $z$, A1 and A3 hold exactly and A2 holds for all directions in both blocks. Yet if $\nu<\omega\sqrt{8\epsilon_0/(1-\epsilon_0)}\approx0.40\,\omega$, then $\hat s=2$ and Step~1 misses $z$ in population; with $\lambda=100$, $\omega=\lambda/3$ and $\nu=10$, the direction $z$ has posterior standard deviation 0.30 and is still missed. Indeed, $\mathrm{He}_2(q_1)\mathrm{He}_2(q_2)$ belongs to the fourth chaos, so $g=0$, $\Hb=\diag(\lambda,\lambda,\nu)$ on $(q_1,q_2,z)$, and $E^{\rm tot}_{\ge3}=4\omega^2$. For $v=(\cos\theta,\sin\theta)$ in the $(q_1,q_2)$ plane, $E_1(v)+E_2(v)=\tfrac12\lambda^2$ by \eqref{eq:E12}. With $t$ the coordinate along $v^\perp=(-\sin\theta,\cos\theta)$, a direct computation gives $\E[\mathrm{He}_2(q_1)\mathrm{He}_2(q_2)\mid t]=\sin^2\theta\cos^2\theta\,\mathrm{He}_4(t)$, so the fourth-order part not involving $v$ has energy $24\omega^2\sin^4\theta\cos^4\theta$ and
\[
E_{\ge3}(v)=E_4(v)=4\omega^2(1-6\sin^4\theta\cos^4\theta)\le4\omega^2\le\tfrac12\lambda^2\qquad(\omega\le\lambda/(2\sqrt2)),
\]
so A2 holds in the N block; $z$ has no term of order three or more, and A3 holds with $\ell=0$, $B=0$, $H=\nu$. With $s_i=q_i^2$, the $(q_1,q_2)$ part of $\Phi$ equals $(\tfrac\lambda2-\omega)(s_1+s_2)+\omega s_1s_2-\lambda+\omega\ge-\lambda+\omega$. In the population greedy path, the first two steps reduce the validation error by at least the fraction $\tfrac12\lambda^2/(2\lambda^2)=\tfrac14$, and the third by $\tfrac12\nu^2/(\tfrac12\nu^2+4\omega^2)$, which is below $\epsilon_0$ if and only if $\nu^2<8\epsilon_0\omega^2/(1-\epsilon_0)$; then $\hat s=2$ and $\cS_{\rm pop}=\mathrm{span}(q_1,q_2)$.

\subsection{Step 2: identifiability, test and basis}\label{app:step2}
Fix the outcome of Step~1, whose data are independent of Step~2. Let $\bar\Phi(w)=\E_\xi[\Phi(F_{\cS}(\mu_{\cS}+A_{\cS}^{-1/2}w)+P_{\cS^\perp}\xi)]$; since $\mathrm{He}_k(w)$ is independent of $\xi$, $\mathcal T_k=\E[\bar\Phi(w)\mathrm{He}_k(w)]$. Let $E^{\rm ref}_k(v)$ be the energies of $\bar\Phi$ under $\Nor(0,I_m)$ (Section~\ref{sec:setup} with $\bar\Phi$ and $w$ in place of $\Phi$ and $x$), $\Pi_2$ the $L^2$ projection onto polynomials of degree at most two, and $\mathcal N^*\subseteq\R^m$ the smallest subspace such that $\bar\Phi-\Pi_2\bar\Phi$ depends only on $P_{\mathcal N^*}w$ (it exists since intersections of such subspaces are such subspaces). Because the two halves are independent, $\E M_{3,4}=\mathcal T_3\mathcal T_3\T/3!+\mathcal T_4\mathcal T_4\T/4!$ up to $O(1/n_2)$ terms from centring the residuals by the full-sample mean; using the same data twice would add a positive semidefinite term proportional to the estimation variance, giving positive energy also to uninformative directions.

\begin{proposition}\label{prop:step2}
Let $\Phi$ be square-integrable under the Step~2 design, so that $\bar\Phi\in L^2(\Nor(0,I_m))$ (this does not follow from $\Phi\in L^2(\Nor(0,I_d))$, as the floor on $A_{\cS}$ widens the design).
\begin{enumerate}[label=(\roman*)]
\item $\range(\E M_{3,4})=\range(\mathcal T_3)+\range(\mathcal T_4)\subseteq\mathcal N^*$, and for unit $v$,
$\tfrac13E^{\rm ref}_3(v)+\tfrac14E^{\rm ref}_4(v)\le v\T(\E M_{3,4})v\le E^{\rm ref}_3(v)+E^{\rm ref}_4(v)$.
Hence $\range(\E M_{3,4})=\mathcal N^*$ if and only if every nonzero direction in $\mathcal N^*$ has $E^{\rm ref}_3+E^{\rm ref}_4>0$.
\item For every sample, after $j$ selections the statistic equals $\lambda_{j+1}(M_{3,4})$, the $(j+1)$th largest eigenvalue ($M_{3,4}$ need not be positive semidefinite), and the selected directions are eigenvectors of $M_{3,4}$. If $\dim\mathcal N^*=r$ and $B_*$ is an orthonormal basis of $\mathcal N^{*\perp}$, then $\lambda_{r+1}(M_{3,4})\le\lambda_{\max}(B_*\T M_{3,4}B_*)$, whose right side involves only $B_*\T\mathcal T_k^{(h)}$, which have mean zero.
\item Suppose the third- and fourth-order coefficient tensors restricted to $\range(E)$ are orthogonally decomposable, $\mathcal C_k=\sum_{i=1}^ra_{ik}o_i^{\otimes k}$ ($k=3,4$) with $\{o_i\}$ an orthonormal basis of $\R^r$ and, for each $i$, $a_{i3}\ne0$ or $a_{i4}\ne0$ (for example, the non-quadratic part of $\bar\Phi$ within $\range(E)$ is a sum of functions of directions orthogonal in $w$). Then the population contrast $\mathcal J(R)=\sum_j[\mathcal C_3(R_j^{\otimes3})^2/3!+\mathcal C_4(R_j^{\otimes4})^2/4!]$ satisfies $\mathcal J(R)\le\sum_i(a_{i3}^2/3!+a_{i4}^2/4!)$, with equality if and only if the columns of $R$ are a permutation of $\{\pm o_i\}$.
\item If $\Phi=\phi_0(a)+(\ell^*+B^*a)\T c+\tfrac12c\T H^*c$ in some coordinates $a=W_N^{*\top}x$, $c=W_Q^{*\top}x$ with orthonormal columns, then $W_Q^{*\top}g=\ell^*$, $W_Q^{*\top}\Hb W_N^*=B^*$, $W_Q^{*\top}\Hb W_Q^*=H^*$, and $\phi_0$ enters only $W_N^{*\top}g$ and $W_N^{*\top}\Hb W_N^*$, without smoothness; if the eigendirections truncated by Step~1 are orthogonal to $\range(W_Q^*)$, the truncated $\Hb_{\hat s}$ gives the same $B^*$ and $H^*$.
\end{enumerate}
\end{proposition}

\begin{proof}
(i) $\bar\Phi-\Pi_2\bar\Phi=\sum_{k\ge3}\bar\Phi_k$ with $\bar\Phi_k=\langle\mathcal T_k,\mathrm{He}_k(w)\rangle/k!$. Conditioning on $P_Vw$ commutes with the chaos projections, so the sum depends only on $P_Vw$ if and only if each $\mathcal T_k\in V^{\otimes k}$, that is $\range(\mathcal T_k)\subseteq V$; thus $\mathcal N^*=\sum_{k\ge3}\range(\mathcal T_k)$, and $\range(XX\T+YY\T)=\range X+\range Y$. In an orthonormal basis whose first vector is $v$, grouping the entries of the symmetric tensor by multi-index $\alpha$ gives $\|\mathcal T_k\times_1v\|^2/k!=\sum_{|\alpha|=k}(\alpha_v/k)\mathcal T_\alpha^2/\alpha!$, while $E^{\rm ref}_k(v)=\sum_{|\alpha|=k,\alpha_v\ge1}\mathcal T_\alpha^2/\alpha!$; as $1\le\alpha_v\le k$, $E^{\rm ref}_k(v)/k\le\|\mathcal T_k\times_1v\|^2/k!\le E^{\rm ref}_k(v)$, and summing over $k=3,4$ gives the bounds. If $\range(\E M_{3,4})\ne\mathcal N^*$, a nonzero $v\in\mathcal N^*$ orthogonal to $\range(\E M_{3,4})$ has $v\T(\E M_{3,4})v=0$, so $E^{\rm ref}_3(v)=E^{\rm ref}_4(v)=0$ by the lower bound; the converse follows from the upper bound.
(ii) If the $j$ selected directions are the leading $j$ eigenvectors of $M_{3,4}$, their complement $\range(B)$ is invariant, the eigenvalues of $B\T M_{3,4}B$ are the remaining ones, the largest is $\lambda_{j+1}(M_{3,4})$, and its eigenvector $v$ satisfies $M_{3,4}Bv=B(B\T M_{3,4}B)v=\lambda_{j+1}Bv$; induction gives the claim. By the Courant--Fischer theorem, $\lambda_{r+1}(M_{3,4})=\min_{\dim V=r}\max_{x\perp V,|x|=1}x\T M_{3,4}x\le\max_{x\perp\mathcal N^*,|x|=1}x\T M_{3,4}x$, and $B_*\T\mathcal T_k=0$ because $\mathcal T_k\in(\mathcal N^*)^{\otimes k}$.
(iii) Let $c_{ji}=o_i\T R_j$, the entries of an orthogonal matrix. By Cauchy--Schwarz and $\sum_ic_{ji}^4\le(\sum_ic_{ji}^2)^2=1$,
\begin{gather*}
\Big(\sum_ia_{i3}c_{ji}^3\Big)^2\le\Big(\sum_ia_{i3}^2c_{ji}^2\Big)\Big(\sum_ic_{ji}^4\Big)\le\sum_ia_{i3}^2c_{ji}^2,\\
\Big(\sum_ia_{i4}c_{ji}^4\Big)^2\le\sum_ia_{i4}^2c_{ji}^4\le\sum_ia_{i4}^2c_{ji}^2 ,
\end{gather*}
and summing over $j$ with $\sum_jc_{ji}^2=1$ gives the bound. If row $j$ is not a signed unit vector, then $\sum_ic_{ji}^4<1$, and equality forces $c_{ji}=0$ whenever $a_{i3}\ne0$ or $a_{i4}\ne0$, that is for all $i$, a contradiction. Distinct energies are not needed.
(iv) A priori $a$, $c$ and the remaining coordinates are independent standard normal. Then $\E[(\ell^*+B^*a)\T c\,c]=\ell^*$, $\E[(\ell^*+B^*a)\T c\,ca\T]=B^*$ and $\E[\tfrac12c\T H^*c(cc\T-I)]=H^*$, whereas $\phi_0(a)$ is independent of $c$ and odd moments of $c$ vanish, so $\phi_0$ does not contribute to $W_Q^{*\top}g$, $W_Q^{*\top}\Hb W_N^*$ or $W_Q^{*\top}\Hb W_Q^*$. If the truncated eigendirections $u_j$ satisfy $W_Q^{*\top}u_j=0$, then $W_Q^{*\top}\Hb_{\hat s}=W_Q^{*\top}\Hb$.
\end{proof}

\paragraph{Asymptotic level (outline).}
Let $r=\dim\mathcal N^*$, $m$ fixed, $n_2\to\infty$ and $\E e^4<\infty$. For a fixed $B$ with $B\T\mathcal T_k=0$ ($k=3,4$), $n_2B\T M_{3,4}B$ converges in distribution to the symmetrised product of two independent Gaussian matrices, and the sign-flipped statistic with the same $B$ converges conditionally to the same limit by the conditional Lindeberg (multiplier) central limit theorem, since flipping preserves the conditional mean zero and the expected conditional covariance; symmetry of $e$ in $w$ is not needed. If $\lambda_r(\E M_{3,4})>0$ and $r\le r_{\max}$, the first $r$ statistics converge to positive constants while the flip quantiles are $O_p(1/n_2)$, so $P(\hat r<r)\to0$. At step $r$, the estimated $B$ converges to a basis of $\mathcal N^{*\perp}$; under regularity conditions such as continuity of the limit distribution at the quantile, the flip quantiles with the estimated and the true basis agree asymptotically, and (ii) gives $\limsup P(\hat r>r)\le\alpha$, the nominal level (about 0.015 here). In finite samples the test is not exact; it can be liberal when the residual has only fifth- or higher-order non-Gaussianity along a direction that grows fast in $|w|$, and the falsely selected direction is then typically that direction.

\subsection{Block structure, positivity and Proposition~\ref{prop:exact}}\label{app:exact}

\paragraph{Proposition~\ref{prop:exact}(i).}
The map $(a,u,\rho)\mapsto x=W_N a+W_Q(\mu(a)+V^{1/2}u)+\rho$ is affine with constant Jacobian, so in these coordinates $\pi$ has density proportional to $\exp\{-\tfrac12|a|^2-\tfrac12|c|^2-\tfrac12|\rho|^2-\Phi(x)\}$, and the (N) proposal, symmetric in $a$ at fixed $(u,\rho)$, has MH ratio $\pi(x')/\pi(x)$, which gives $\alpha_N$. If $\Phi=\phi_0(a)+(\ell+Ba)\T c+\tfrac12c\T Hc$, then with $V^{-1}=I_k+H$ and $\mu(a)=-V(\ell+Ba)$,
\[
\tfrac12|c|^2+(\ell+Ba)\T c+\tfrac12c\T Hc=\tfrac12\big(c-\mu(a)\big)\T V^{-1}\big(c-\mu(a)\big)-\tfrac12\mu(a)\T V^{-1}\mu(a),
\]
and the first term is $\tfrac12|u|^2$. The density factorises as $m(a)e^{-|u|^2/2}e^{-|\rho|^2/2}$, the (N) ratio becomes $m(a')/m(a)$, and $\alpha_U\equiv1$ because $\Phi$ does not depend on $\rho$.

\paragraph{Block structure.}
In (N), $\mu(a')-\mu(a)=-VB(a'-a)=c'-c$, so $u'=V^{-1/2}(c'-\mu(a'))=u$ and $\rho$ is unchanged; in (Q), $a$ and $\rho$ are unchanged and $u'-u=V^{-1/2}S\xi$; in (U), $a$ and $c$, hence $u$, are unchanged. The three acceptance probabilities are Metropolis(--Hastings) ratios for $p(\cdot\mid u,\rho)$, $p(\cdot\mid a,\rho)$ and $p(\cdot\mid a,u)$, where $p$ is the density of $\pi$ in $(a,u,\rho)$. This holds for any $\Phi$.

\paragraph{Acceptance probability of (U).}\label{app:U}
A priori $\rho$ is independent of $(a,c)$ with distribution $\Nor(0,P_U)$. The refresh proposal is this distribution, and the pCN proposal makes $(\rho,\rho')$ an exchangeable Gaussian pair when $\rho\sim\Nor(0,P_U)$; both are reversible with respect to $\Nor(0,P_U)$, so the MH ratio reduces to $e^{\Phi(x)-\Phi(x')}$, and an equal mixture of the two is reversible.

\paragraph{Positivity.}
Let a proposal kernel $Q$ be reversible with respect to a $\sigma$-finite measure $\nu$ and positive on $L^2(\nu)$, let the target be $f\nu$, and let $P$ be the Metropolis kernel with acceptance probability $\min\{1,f(y)/f(x)\}$. For $g\in L^2(f\nu)$,
\[
\langle g,Pg\rangle_{f\nu}=\iint g(x)g(y)\min\{f(x),f(y)\}\,\nu(dx)Q(x,dy)+\int r\,g^2f\,d\nu ,
\]
with $r\ge0$ the rejection probability. Since $\min\{f(x),f(y)\}=\int_0^\infty\mathbf 1\{f(x)>t\}\mathbf 1\{f(y)>t\}\,dt$, the first term equals $\int_0^\infty\langle g\mathbf 1_{A_t},Q(g\mathbf 1_{A_t})\rangle_\nu\,dt\ge0$ with $A_t=\{f>t\}$; here $\int_{A_t}g^2d\nu\le t^{-1}\|g\|^2_{f\nu}$, and Fubini's theorem applies because $|g(x)g(y)|\min\{f(x),f(y)\}\le\tfrac12(g^2f)(x)+\tfrac12(g^2f)(y)$ and $Q$ is $\nu$-reversible. This is the argument of \citet{rudolf2013}, written there for the uniform measure on bounded sets. For (N), fix $(u,\rho)$ and take Lebesgue measure on $a$; the proposal is convolution with the symmetric measure $r^{-1}\sum_j\Nor(0,\sigma_j^2)e_j$, whose Fourier transform $r^{-1}\sum_je^{-\sigma_j^2\omega_j^2/2}$ is nonnegative. (Q) is similar. For (U), take $\nu=\Nor(0,P_U)$: the refresh is the projection onto constants and pCN is the Mehler operator, equal to $(1-b_U^2)^{n/2}\ge0$ on the $n$th chaos. Integrating over the other blocks gives $\langle g,P_bg\rangle_\pi\ge0$ for each block, so $P=\sum_bp_bP_b$ is positive and its spectrum lies in $[0,1]$.

\paragraph{Proposition~\ref{prop:exact}(ii)--(iii).}
By (i) and the form of the acceptance ratios, the block kernels have the stated tensor form; they commute and are self-adjoint. For a nonempty $T\subseteq\{N,Q,U\}$ let $\mathcal H_T$ be the tensor product of $L^2_0$ for the blocks in $T$ and the constants for the others. Then $L^2_0(\pi)=\bigoplus_T\mathcal H_T$, each $\mathcal H_T$ is invariant under all block kernels, and on $\mathcal H_T$, $I-P=\sum_{b\in T}p_b(I-K_b)\succeq(\sum_{b\in T}p_b\gamma_b)I$. The minimum over $T$ is attained by single blocks, and single-block functions give equality, so $\gap(P)=\min_bp_b\gamma_b$; positivity is not needed for this, only for identifying the right and absolute gaps. For a function $f$ of one block $b$ only, with selection probability $p_b$, $P=(1-p_b)I+p_bK_b$ on $\mathcal H_{\{b\}}$, and integrating against the spectral measure $\mu_f$ of $f$ for $K_b$ gives
\begin{align*}
\IACT_P(f)&=\int\frac{2-p_b(1-\lambda)}{p_b(1-\lambda)}\,\frac{\mu_f(d\lambda)}{\Var f}=\frac1{p_b}\int\frac{2}{1-\lambda}\,\frac{\mu_f(d\lambda)}{\Var f}-1\\
&=\frac{\IACT_{K_b}(f)+1}{p_b}-1 ,
\end{align*}
and for sums of single-block functions the asymptotic variances add, because the $\mathcal H_{\{b\}}$ are orthogonal and invariant. For (iii), $K_U$ has eigenvalues $\tfrac12s^n$, $n\ge1$, on $L^2_0$, with $s=\sqrt{1-b_U^2}$, so $\gamma_U=1-\tfrac12s$. As $\Nor(0,I_k)$ is rotation invariant, $\gamma_Q$ depends only on the eigenvalues $s_1\le\dots\le s_k$ of $\Sigma_u$. In the eigenbasis take $f=u_1$ and the proposal $u'=u+\Sigma_u^{1/2}\xi$; then $\mathcal E(f)=\tfrac12\E[\alpha\,s_1\xi_1^2]$. Using $\min\{1,e^\Delta\}\le e^{\Delta/2}$ with $\Delta=\sum_i\Delta_i$, $\Delta_i=-\tfrac12[(u_i+\sqrt{s_i}\xi_i)^2-u_i^2]$, and $\E[e^{\Delta_i/2}\mid\xi_i]=e^{-s_i\xi_i^2/8}$,
\begin{align*}
\gamma_Q\le\mathcal E(f)&\le\tfrac12s_1\big(1+\tfrac{s_1}4\big)^{-3/2}\prod_{i\ge2}\big(1+\tfrac{s_i}4\big)^{-1/2}\le\tfrac12s_1\big(1+\tfrac{s_1}4\big)^{-(k+2)/2}\\
&\le\frac4k\Big(1+\frac2k\Big)^{-(k+2)/2}\le\frac4{ek},
\end{align*}
where the maximum over $s_1$ is at $s_1=8/k$ and the last step uses $(1+t)\log(1+t)\ge t$.

\paragraph{Bounded perturbations.}
If $\zeta=\Phi-\Phi_E$ is bounded for some exact-class $\Phi_E$ with the sampler's $\ell,B,H$, then for the same frozen parameters $\gap(P_\pi)\ge e^{-\mathrm{osc}\,\zeta}\gap(P_{\pi_E})$, $\mathrm{osc}\,\zeta=\sup\zeta-\inf\zeta$. To see this, let $\pi=h\pi_E$ with $h\propto e^{-\zeta}$, so $\inf h/\sup h=e^{-\mathrm{osc}\,\zeta}$. By the representation above, the off-diagonal part of each block kernel is $\min\{f(x),f(y)\}\nu_b(dx)Q_b(x,dy)$, and $f=hf_E$ gives $\min\{f(x),f(y)\}\ge(\inf h)\min\{f_E(x),f_E(y)\}$; as Dirichlet forms involve only off-diagonal parts, $\mathcal E_{P_\pi}(g)\ge(\inf h)\mathcal E_{P_{\pi_E}}(g)$. Since also $\Var_\pi g\le\E_\pi(g-\E_{\pi_E}g)^2\le(\sup h)\Var_{\pi_E}g$, $\gap(P_\pi)\ge e^{-\mathrm{osc}\,\zeta}\gap(P_{\pi_E})$.

\subsection{Invariance of (Q) and diminishing adaptation}\label{app:Q}
Given $(a,\rho)$, the target of $c$ is proportional to $\exp\{-\tfrac12|c|^2-\Phi(W_N a+W_Q c+\rho)\}$, and for fixed $S$ the proposal $c'=c+S\xi$ is symmetric, so each $P_{Q,S}$ leaves $\pi$ invariant. The RAM update is $S_{n-1}^{-1}\Sigma_nS_{n-1}^{-\top}=I_k+\varepsilon_n\tilde\xi\tilde\xi\T$ with $\tilde\xi=\xi/|\xi|$ and $\varepsilon_n=\eta_n(\alpha_n-\alpha^*)\in[-\eta_n\alpha^*,\eta_n(1-\alpha^*)]$; as $\eta_n\le1$, $1+\varepsilon_n\ge1-\alpha^*>0$ and $\Sigma_n$ stays positive definite. Total variation between Gaussians is invariant under invertible linear maps, so consecutive proposals differ by $\|\Nor(0,1)-\Nor(0,1+\varepsilon_n)\|_{\rm TV}=O(|\varepsilon_n|)$. Two Metropolis kernels with the same target and symmetric proposals $q_1,q_2$ satisfy $\|P_1(x,\cdot)-P_2(x,\cdot)\|_{\rm TV}\le\|q_1(x,\cdot)-q_2(x,\cdot)\|_{\rm TV}$, because the acceptance probability does not depend on the proposal and lies in $[0,1]$. Weighting by the selection probability $p_Q$ of (Q) gives $\sup_x\|P_{t+1}(x,\cdot)-P_t(x,\cdot)\|_{\rm TV}\le p_Q\|\Nor(0,1)-\Nor(0,1+\varepsilon_n)\|_{\rm TV}\to0$, since $|\varepsilon_n|\le\eta_n\max\{\alpha^*,1-\alpha^*\}\to0$; here $p_Q=0.3$ (0.5 if $r=0$; 0.27 or 0.45 with a global scale, the supremum then also over $\tau$). With a global scale, the proposals $c'=c+V_\tau^{1/2}S\xi$ differ by one more invertible map and the bound is unchanged up to the weight.

In the exact class more can be said about the adaptation. With $\Sigma_n=S_nS_n\T$ and $\Delta=c'-c$, the RAM update reads $\Sigma_n=\Sigma_{n-1}+\eta_n(\alpha_n-\alpha^*)\Delta\Delta\T/(\Delta\T\Sigma_{n-1}^{-1}\Delta)$; it does not depend on the choice of square root and is equivariant under invertible linear maps. In the exact class the (Q) acceptance depends only on $u$ and $u'$, and (N) and (U) do not change $u$. Hence, counted in (Q) updates, $(u_n,V^{-1/2}\Sigma_nV^{-\top/2})$ has the law of RAM on $\Nor(0,I_k)$ started from $(2.38^2/k)I_k$, independently of $\phi_0,\ell,B,H$, $d$ and the evolution of $(a,\rho)$. Consequently $(a,\rho)$ is itself a time-homogeneous reversible Markov chain with invariant distribution $\pi_a\otimes\Nor(0,P_U)$, so standard Markov chain theory applies to ergodic averages of $(a,\rho)$, which include the N and U directions, whatever the adaptation does. Functions of $u$ reduce to the stability of RAM on $\Nor(0,I_k)$, Vihola's assumption, which to our knowledge is unproved for $k\ge2$; on the event of stability, $S_nS_n\T\to\theta_k^*V$ with $\theta_k^*$ the constant giving mean acceptance $\alpha^*$ on $\Nor(0,I_k)$. With a global scale, (H) redraws $u$ from $\Nor(0,I_k)$ in the exact class (Proposition~\ref{prop:tau}), so the kernel of $u$ is bounded below by $p_H\Nor(0,I_k)$ for every $S$; with diminishing adaptation this gives ergodicity and a weak law of large numbers for bounded functions of $u$ without the stability of RAM \citep{roberts2007}. None of this applies outside the exact class, where the (Q) acceptance depends on $(a,\rho)$.

\subsection{Comparison with pCN}\label{app:pcn}
Let $\Phi$ depend on $x$ only through $P_{S'}x$ for a subspace $S'$, let $v_0\perp S'$ be a unit vector, let $\vartheta=\sqrt{1-\beta^2}$ for the pCN parameter $\beta$, and let $\bar\alpha(\beta)$ be the stationary mean acceptance rate. Under $\pi$ and the proposal, $(v_0\T x,v_0\T x')$ is standard bivariate normal with correlation $\vartheta$ and independent of the acceptance probability, so $f=\mathrm{He}_n(v_0\T x)$ has Dirichlet form $\bar\alpha(\beta)(1-\vartheta^n)\Var f$ and
\begin{gather*}
\gap(\mathrm{pCN}_\beta)\le\bar\alpha(\beta)\big(1-\sqrt{1-\beta^2}\big)\le\bar\alpha(\beta)\,\beta^2,\\
\IACT\big(\mathrm{He}_n(v_0\T x)\big)\ge\frac{2}{\bar\alpha(\beta)(1-\vartheta^n)}-1 .
\end{gather*}
Moreover $\bar\alpha(\beta)$ is at most the Bhattacharyya coefficient between $\pi(dx)q(x,dy)$ and its reverse, which factorises if the posterior factorises along orthogonal directions, with factor at most $\{1+\beta^2\lambda^2/(4(1+\lambda))\}^{-1/2}$ for a direction with likelihood $\tfrac12\lambda(v\T x-m)^2$ and at most 1 otherwise. Hence, if the posterior factorises along orthogonal directions and at least $n_\Lambda\ge3$ of them have such a likelihood with $\lambda\ge\Lambda\ge1$, then $\gap(\mathrm{pCN}_\beta)\le16/\{e(n_\Lambda-2)\Lambda\}$ for every $\beta$: the stiffest directions force $\beta$ to be small and all weakly informed and uninformative directions slow down, whereas in the exact class with $SS\T\propto V$ the Q and U gaps do not depend on the stiffness. For $\Phi_G$ of Section~\ref{sec:numerics} under $\Nor(0,I_{200})$, the ten Gaussian directions alone give $\gap(\mathrm{pCN}_\beta)\le1.25\times10^{-3}$ for every $\beta$.

\subsection{Further results on the global scale, and proofs}\label{app:tauproofs}

Under the assumptions of Proposition~\ref{prop:tau}, two further properties hold.
\begin{enumerate}[label=(\roman*),start=3]
\item The kernel $P_S$ of the whole chain satisfies $\min\{\gap(M),p_H\}\le\gap(P_S)\le\gap(M)$, and $P_S^nf=M^nf$ for functions $f$ of $(a,\tau)$.
\item If $r=0$, then $\gap(P_S)=p_H\gap(\mathcal K_H)$ and, for functions $f$ of $\tau$, the IACT per call is $(\IACT_{\mathcal K_H}(f)+1)/p_H-1$. If moreover $p(\eta\mid y)$ is log-concave and $\mathcal K_H$ is ideal slice sampling, then $\gap(\mathcal K_H)\ge\tfrac12$, so $\gap(P_S)\ge p_H/2=0.05$ and every $f\in L^2(\pi)$ has IACT per call at most 39, whatever $d$, $k$, the stiffness of the Q directions and the data.
\end{enumerate}
The bound $\tfrac12$ in (iv) is due to \citet{natarovskii2021}. The implemented stepping-out and shrinkage coincide with ideal slice sampling on unimodal targets whenever both ends step out of the slice before the budget is exhausted.

\begin{proposition}\label{prop:centered}
Let $\pi(x,\tau)$ be as in Section~\ref{sec:tau} with arbitrary $\Phi$ and $\E_\pi\tau^2<\infty$. Consider a random-scan mixture of $\pi$-reversible kernels in which the components that change $\tau$ update it given all of $x$ (total selection probability $p_\tau$; each reversible with respect to $\pi(\tau\mid x)$, e.g.\ exact Gibbs, Metropolis or slice sampling) and the other components do not change $\tau$. Then
\[
\gap\ \le\ 2p_\tau\,\frac{\E_\pi[\Var(\tau\mid x)]}{\Var_\pi(\tau)},
\]
and the factor 2 can be dropped if the $\tau$ kernels are positive operators (e.g.\ exact Gibbs or ideal slice sampling). For an $\mathrm{InvGamma}(\alpha_\tau,\beta_\tau)$ hyperprior, $\Var(\tau\mid x)=\E[\tau\mid x]^2/(\alpha_\tau+d/2-2)$ and
\[
\gap\ \le\ \frac{4p_\tau(1+\mathrm{CV}^{-2})}{d+2\alpha_\tau-4},\qquad \mathrm{CV}=\frac{\mathrm{sd}_\pi(\tau)}{\E_\pi\tau}.
\]
\end{proposition}

Proposition~\ref{prop:centered} holds for any $\Phi$. If the number of informed directions is fixed, the posterior of $\tau$ and hence CV do not depend on $d$, and the bound decreases as $O(1/d)$; under an inverse-gamma hyperprior the coefficient of variation of $\tau$ given $x$ is exactly $(\alpha_\tau+d/2-2)^{-1/2}$, the $\sqrt{2/d}$ of Section~\ref{sec:tau}. \citet{agapiou2014} obtained the same order asymptotically for centred Gibbs sampling in linear Gaussian inverse problems; Proposition~\ref{prop:centered} is non-asymptotic and allows nonlinear $\Phi$. The comparison concerns the order in $d$; at fixed $d$ it does not say which update is faster.

\paragraph{Proof of Proposition~\ref{prop:tau} and of (iii)--(iv).}
(i) For fixed $\tau$, $(a,u,\tilde\rho)\mapsto x=W_N a+W_Q(\mu_\tau(a)+V_\tau^{1/2}u)+\sqrt\tau\,\tilde\rho$ is block-triangular and affine with Jacobian $|V_\tau|^{1/2}\tau^{d_U/2}$, $d_U=d-r-k$, independent of $a$. Completing the square, $|c|^2/(2\tau)+b(a)\T c+\tfrac12c\T H_+c=\tfrac12|u|^2-\tfrac12b(a)\T V_\tau b(a)$, and $|\rho|^2/(2\tau)=\tfrac12|\tilde\rho|^2$. Since $\tau^{-d/2}\tau^{d_U/2}|V_\tau|^{1/2}=\tau^{-r/2}(|V_\tau|/\tau^k)^{1/2}$,
\[
\pi(a,\tau,u,\tilde\rho)\ \propto\ p(\tau)\,\tau^{-r/2}e^{-|a|^2/(2\tau)-\phi_0(a)}\,Z_Q(a;\tau)\ e^{-|u|^2/2}\,e^{-|\tilde\rho|^2/2},
\]
so $\pi(\tau\mid a)\propto p(\tau)\Nor(a;0,\tau I_r)Z_Q(a;\tau)$, which is $\hat p$ in $\eta$, and $c,\rho$ given $(a,\tau)$ are as stated. As $R=\phi_0(a)$, $\alpha_H\equiv1$.
(ii) In $(a,\tau,u,\tilde\rho)$, (N) fixes $(u,\rho,\tau)$, hence $\tilde\rho$, with acceptance ratio $\pi_{a\tau}(a',\tau)/\pi_{a\tau}(a,\tau)$, where $\pi_{a\tau}$ is the first factor above; (Q) moves $u$ by $S\xi$ with acceptance depending on $u$ only; (U) moves $\tilde\rho$ and is always accepted; (H) draws $\tau'\sim\mathcal K_H(\tau,\cdot)$ and fresh $u'$, $\tilde\rho'$. So $(a,\tau)$ evolves by $M$, and $\pi_{a\tau}$ does not involve $d$.
(iii) Let $\mathcal H_1$ be the functions of $(a,\tau)$ and $\mathcal H_\perp=\{f:\E_\pi[f\mid a,\tau]=0\}$. All four kernels preserve both subspaces; $P_S=M$ on $\mathcal H_1$; on $\mathcal H_\perp$ the output of (H) depends only on $(a,\tau')$, so $P_H=0$ there and $\|P_S|_{\mathcal H_\perp}\|\le1-p_H$. For $f=f_1+f_2\in L^2_0(\pi)$, $\mathcal E_{P_S}(f)=\mathcal E_M(f_1)+\langle f_2,(I-P_S)f_2\rangle\ge\gap(M)\|f_1\|^2+p_H\|f_2\|^2$; functions in $\mathcal H_1$ give the upper bound and $P_S^nf=M^nf$.
(iv) For $r=0$, $M=(1-p_H)I+p_H\mathcal K_H$, so $\gap(M)=p_H\gap(\mathcal K_H)\le p_H$, since the gap of a reversible kernel is at most 1 when $\pi(\tau\mid y)$ has no atoms (take $f=\mathbf 1_A-\pi(A)$ with $\pi(A)\to0$), and (iii) gives $\gap(P_S)=p_H\gap(\mathcal K_H)$; the IACT formula follows as in Appendix~\ref{app:exact}. The spectral gap of ideal slice sampling on a one-dimensional density depends only on the lengths $v(t)=|\{\varrho>t\}|$ of its level sets \citep{natarovskii2021}. For a log-concave $\varrho$ these are intervals and $v$ is concave in $-\log t$, so the symmetric decreasing rearrangement of $\varrho$ is log-concave with the same $v$; for symmetric one-dimensional log-concave densities, \citet{natarovskii2021} show that the gap is at least $\tfrac12$. Hence $\gap(P_S)\ge p_H/2$ and $\IACT(f)\le2/\gap-1\le39$.

\paragraph{Proof of Proposition~\ref{prop:centered}.}
Take $f(x,\tau)=\tau$. Components that do not change $\tau$ do not contribute to $\mathcal E(f)$. For a $\mu$-reversible $L$, $\mathcal E_L(g)=\Var_\mu g-\mathrm{Cov}(g(X),g(X'))\le2\Var_\mu g$, and $\mathcal E_L(g)\le\Var_\mu g$ if $L$ is positive. Hence $\mathcal E(f)\le2p_\tau\E_\pi[\Var(\tau\mid x)]$, and $\gap\le\mathcal E(f)/\Var_\pi f$. Since $\Phi$ does not involve $\tau$, $\tau\mid x\sim\mathrm{InvGamma}(\alpha_\tau+d/2,\beta_\tau+|x|^2/2)$ and $\Var(\tau\mid x)=\E[\tau\mid x]^2/(\alpha_\tau+d/2-2)$; moreover $\E[\E[\tau\mid x]^2]\le\E\tau^2=\Var(\tau)(1+\mathrm{CV}^{-2})$.

\section{Hermite energies of the examples}\label{app:energy}
If $\Phi=\sum_i\phi_i(v_i\T x)$ for orthonormal $v_i$, then $E_k(v_i)=\E[\phi_i(z)\mathrm{He}_k(z)]^2/k!$ with $z\sim\Nor(0,1)$. Each $\Phi$ is a sum of terms depending on one coordinate or on a pair $(z_1,q_j)$. For each true direction $i$, the main effect $F_i(t)=\E[\Phi\mid z_i=t]$ is available in closed form; its Hermite energies up to order 200 are computed on a grid of $2.4\times10^6$ points on $[-12,12]$, energies of order three or more as the variance minus the first two. Mixed terms add $(0.3\lambda_j)^2$ (second order, M2) or $2(b_j\lambda_j)^2$ (third order, $-b_j\lambda_j\mathrm{He}_2(z_1)\mathrm{He}_1(q_j)$, M3) to each direction involved.

A2 (Section~\ref{sec:setup}) is not imposed on mixtures of N and Q directions, whose first- and second-order energies can cancel. Within the plane of N1 and N2, where the energies follow from the Hermite coefficient tensors of the terms of $\Phi$, the maximum of $E_{\ge3}/(E_1+E_2)$ is 1.02, at N2, in all examples.

Table~\ref{tab:energy} gives the energies of the N directions and of the coupled directions of M3 under the learning design; for M4 the design uses $\tau=1$. The Gaussian directions of M1, M2, M4 and all directions of Heat have no energy of order three or more. A3 holds exactly in the true coordinates of M1 and M2 (in M2 the conditional means are linear in $z_1$) and, given $\tau$, of M4 and Heat. In M3 the conditional means of Q7--Q10 are quadratic in $z_1$: replacing them by their best affine approximation in $z_1$ under the posterior costs an expected Kullback--Leibler divergence of 0.020 nats, and 36\% of the posterior variance of Q7 comes from the variation of its conditional mean.

\begin{table}[htbp]
\centering\small\setlength{\tabcolsep}{4pt}
\caption{Hermite energies under the learning design $\Nor(0,I_{200})$ and the ratio $E_{\ge3}/(E_1+E_2)$.}\label{tab:energy}
\begin{tabular}{llrrrr}
\toprule
Example & Direction & $E_1$ & $E_2$ & $E_{\ge3}$ & Ratio \\
\midrule
M1, M4 & N1 (bimodal) & 0 & 1.85 & 0.995 & 0.54 \\
M2 & N1 (bimodal, linear coupling) & $1.20\times10^4$ & $1.22\times10^5$ & 0.995 & $7\times10^{-6}$ \\
M3 & N1 (unimodal, curved coupling) & 48.8 & 16.0 & 14.6 & 0.22 \\
M3 & Q7, Q9 (curved coupling) & 3.84, 0.346 & 50, 4.5 & 8, 0.72 & 0.149 \\
M3 & Q8, Q10 (curved coupling) & 16.2, 1.45 & 12.5, 1.125 & 2, 0.18 & 0.070 \\
M1--M4 & N2 (spike) & 0 & 1.77 & 1.81 & 1.02 \\
\bottomrule
\end{tabular}
\end{table}

\section{Implementation}\label{app:impl}

\begin{table}[htbp]
\centering\small
\caption{Settings of the NQU sampler (all runs).}\label{tab:impl}
\begin{tabularx}{\textwidth}{>{\raggedright\arraybackslash}p{2.4cm}>{\raggedright\arraybackslash}X}
\toprule
Item & Value \\
\midrule
Step 1 & $n_1=30\,000$ (also $20\,000$); 20\% validation; ScaledGD with Armijo backtracking, at most 500 iterations per rank, relative tolerance $10^{-9}$; maximal rank 40; stop if the validation error falls by less than 2\%; $\delta=0.1$ \\
Step 2 & $n_2=20\,000$ (also $30\,000$); floor 0.5 on the eigenvalues of $A_{\cS}$; halves in order; 200 sign flips, 0.99 quantile with linear interpolation (level about 0.015); $r_{\max}=5$; $W_N=\orth(F_{\cS}A_{\cS}^{-1/2}E)$ \\
Basis within N & For $r\ge2$: cross-fitted $\mathcal J(R)$ maximised by Jacobi rotations (grid of 180 angles in $[0,\pi/2)$, then bounded one-dimensional refinement; rotations accepted only if $\mathcal J$ increases; stop after a sweep with relative gain $\le10^{-10}$ or 50 sweeps); $W_N\gets W_N\polar(W_N\T\Gamma)$ \\
Step 3 & Eigenvalues of $I_k+H$ restricted to $[0.5,10^6]$; with a global scale, negative eigenvalues of $H$ set to 0 \\
Sampling & $T=450\,000$, warm-up 45\,000; selection $(0.4,0.3,0.3)$, $(0,0.5,0.5)$ if $r=0$; with a global scale $(0.36,0.27,0.27,0.10)$, $(0,0.45,0.45,0.10)$ if $r=0$; start $F_{\cS}\mu_{\cS}+P_{\cS^\perp}\xi$, $\tau_0=1$ \\
(N) & $\sigma_j$ initially 1; target acceptance 0.35; epochs of 100 calls growing by the factor 1.15; adjusted only after epochs with at least 20 (N) proposals, for coordinates with at least 10 \\
(Q) & RAM with Gaussian proposals; $S_0S_0\T=(2.38^2/k)V$; $\eta_n=\min\{1,kn^{-2/3}\}$; $\alpha^*=0.234$; adaptation throughout; with a global scale in $u$ with $S_0=(2.38/\sqrt k)I_k$ \\
(U) & Equal mixture of refresh and pCN; $b_U$ initially 1 and kept at 1 while refreshes are accepted at rate at least 0.3, otherwise adjusted every 100 counted updates towards acceptance 0.3, $b_U\in[10^{-3},1]$ \\
(H) & Slice sampling in $\log\tau$ (stepping out with width 1, at most 50 steps, shrinkage) on $\hat p(\eta\mid a)$; redraw $c$ and $\rho$; accept with $\min\{1,e^{R(x)-R(x')}\}$ \\
\bottomrule
\end{tabularx}
\end{table}

\section{Experimental details}\label{app:exp}

Table~\ref{tab:baselines} gives the settings of the baselines. The threshold 0.01 of AS-MCMC and QPHD equals $\delta^2$ of the NQU sampler and roughly matches the DILI threshold 0.1, as the matrices here scale with squared curvature. The original AS-MCMC \citep{constantine2016} draws the inactive variables from the prior and is approximate; we sample them by MCMC. Dense NUTS used $7.06$--$8.76\times10^5$ calls in warm-up alone, more than the budget, so its ESS within the budget is zero; Table~\ref{tab:main} ignores the budget for it, counting its warm-up and estimating its rate from the 2000 iterations after warm-up. For M4 and Heat, pCN and NUTS run on the non-centred parameter $x=\tau^{1/2}\xi$ (pCN on $\xi$ with probability 0.9, random walk on $\log\tau$ otherwise), and AS-MCMC and QPHD use the (H) update ($r=0$). DILI, gpCN, SAFES and SAFES-P, which assume a Gaussian prior, are applied to the non-centred parameter $(\xi,w)\sim\Nor(0,I_{201})$ with $x=\tau^{1/2}\xi$ and $\tau=F^{-1}(\Psi(w))$, where $F$ is the $\mathrm{InvGamma}(3,2)$ distribution function and $\Psi$ the standard normal one; DILI uses Gauss--Newton Hessians for Heat, whose misfit is of least-squares form, and exact Hessians for M4. On the posterior axes $V_{\rm full}\T x$, the leading columns of $V_{\rm full}$ are the informed directions of the example (for Heat, the right singular vectors of $G$) completed by a fixed orthonormal basis; for M4 and Heat both coordinate systems include $\tau$, giving 402 quantities instead of 400. The baselines use their published or default settings rather than tuning for these examples, and the term $10^6$ divided by the minimum ESS per 1000 calls extrapolates the rate observed in the ESS window. Table~\ref{tab:bias} counts the chains flagged as biased.

\begin{table}[htbp]
\centering\small
\caption{Baselines. The last column gives the calls before the ESS window (NQU sampler: 95\,000, i.e.\ learning 50\,000 and warm-up 45\,000).}\label{tab:baselines}
\begin{tabularx}{\textwidth}{>{\raggedright\arraybackslash}p{1.8cm}>{\raggedright\arraybackslash}p{1.9cm}>{\raggedright\arraybackslash}X>{\raggedright\arraybackslash}p{1.9cm}}
\toprule
Method & Derivatives & Settings & Before window \\
\midrule
pCN & none & Full-space pCN \citep{cotter2013}; step tuned to acceptance 0.25 over the first $5\times10^4$ calls, then fixed & 50\,000 \\
gpCN & none & \citet{rudolf2018} with the empirical mean and covariance of the chain, adapted over $2.5\times10^5$ calls (every 5000 calls from call 20\,000), step tuned to acceptance 0.25, then fixed; first 10\% of sampling discarded & 275\,000 \\
QPHD & none & Full quadratic model $c+b\T x+\tfrac12x\T Ax$ (20\,301 coefficients) fitted by least squares on $5\times10^4$ prior points; $W_Q$ spanned by eigenvectors of $bb\T+A^2$ with eigenvalues $\ge0.01$ \citep{constantine2016joss,constantine2017}; Gaussian from the quadratic model on $W_Q$; sampling as AS-MCMC & 95\,000 \\
SAFES, SAFES-P & none & Algorithms 3.3 and 3.4 of \citet{dunlop2022}. From the paper: 40 particles started from the prior, $\lambda=0.2$ fixed, $\beta$ tuned to keep the acceptance rate in $(0.15,0.3)$, first 25\% of sweeps burn-in, $M=20$ for SAFES-P; particles updated one after another within a sweep. Our choices, not specified in the paper: $\beta$ starts at 0.5 and every 10 sweeps is multiplied by 0.8 (acceptance $<0.15$) or 1.25 ($>0.3$, capped at 1), during burn-in only, with no lower bound. ESS from the pooled autocovariance of the particle paths & 125\,000 \\
DILI & gradient, Hessian & \citet{cui2016}: LIS from the posterior-averaged Hessian (eigenvalues $\ge0.1$); operator-weighted MGLI-Langevin proposal, time step tuned to acceptance 0.574; $5\times10^4$ adaptive calls (LIS rebuilt every 5000 calls from Hessians at every 100th state) plus 451 Hessians, then frozen; first 10\% of sampling discarded & 95\,405 \\
NUTS (diagonal, dense) & gradient & NumPyro 0.22 \citep{phan2019}, float64, Stan-style adaptation, target acceptance 0.8, maximum tree depth 10, 1000 warm-up iterations; diagonal or dense mass matrix; started from the prior & diagonal $0.70$--$1.54\times10^5$; dense $7.06$--$8.76\times10^5$ \\
AS-MCMC (exact) & gradient & Gradients at $5\times10^4$ prior points (the first $3\times10^4$ shared with Step~1); $W_Q$ spanned by eigenvectors of $\overline{\nabla\Phi\nabla\Phi\T}$ with eigenvalues $\ge0.01$ \citep{constantine2016}; Gaussian model on $W_Q$ from Stein's identity; sampling stage of the NQU sampler with $r=0$ & 95\,000 \\
\bottomrule
\end{tabularx}
\end{table}

\begin{table}[htbp]
\centering\small
\caption{Number of biased chains ($\max|z|>4$) out of 8, posterior axes / original coordinates.}\label{tab:bias}
\begin{tabular}{llllll}
\toprule
Method & M1 & M2 & M3 & M4 & Heat \\
\midrule
NQU (30k/20k) & 1 / 0 & 0 / 0 & 0 / 1 & 0 / 1 & 0 / 0 \\
NQU (20k/30k) & 0 / 1 & 1 / 0 & 0 / 0 & 0 / 0 & 0 / 0 \\
DILI & 1 / 0 & 1 / 1 & 0 / 1 & 1 / 0 & 0 / 1 \\
NUTS (diagonal) & 1 / 2 & 1 / 2 & 1 / 1 & 2 / 2 & 0 / 0 \\
NUTS (dense) & 3 / 3 & 2 / 2 & 1 / 0 & 0 / 0 & 0 / 0 \\
AS-MCMC (exact) & 1 / 4 & 1 / 0 & 1 / 2 & 0 / 0 & 8 / 8 \\
pCN & 1 / 1 & 5 / 6 & 5 / 2 & 3 / 3 & 8 / 8 \\
gpCN & 8 / 8 & 8 / 8 & 8 / 8 & 8 / 8 & 8 / 8 \\
QPHD & 6 / 7 & 8 / 7 & 8 / 8 & 2 / 3 & 0 / 0 \\
SAFES & 8 / 8 & 8 / 8 & 8 / 8 & 8 / 8 & 8 / 8 \\
SAFES-P & 8 / 8 & 8 / 8 & 8 / 8 & 8 / 8 & 8 / 8 \\
\bottomrule
\end{tabular}
\end{table}

The learning diagnostics are as follows. In finite samples each new column of Step~1 also costs the estimation of about $d-s$ parameters, which raises the relative stopping threshold to about $\epsilon_0+\kappa(d-s)/n_t$ with $n_t=0.8\,n_1$; a fixed-$d$ expansion gives $\kappa\approx1$, and we observed $\kappa$ between 0.8 and 2.4 on synthetic examples. This is an approximation, not a theorem; for $d=200$ and $n_t=24\,000$ it gives about 3\%. The step of Step~1 that adds N2 reduces the population validation error by 18--31\% (M1, M2, M4) or 9.2\% (M3), well above this threshold. Extending the curved coupling of M3 to the six Gaussian directions with $\lambda_j\le50$ ($b_j=0.1$) keeps A2 for N1 and all Gaussian directions but raises the total higher-order energy from 16.4 to 72.4 and lowers this reduction to 2.35\%, above $\epsilon_0$ but below the finite-sample threshold; with $n_1=30\,000$, 7 of 8 seeds leave N2 in U, and with $n_1=120\,000$, 5 of 8 find it; this variant was not sampled. In M1, $\Hb$ has rank 192, the population rank of Step~1 is $\hat s=25$ and its population output has dimension 26 (the extra direction is the part of $g$ in the truncated weakly informed directions); the learned dimensions 23--25 estimate this target (Proposition~\ref{prop:step1}). In M2, about 0.59 of N1 lies in N in population: when N1 moves by one unit, the conditional mean of $q_j$ moves by $0.3\lambda_j/(1+\lambda_j)$. The extra dimension of $\cS$ in M2--M4 is the estimation error of $\hat g$: for homoscedastic residuals of variance $\sigma^2$, the part of $\hat g$ outside $\range(\hat H)$ gives an eigenvalue of about $(d-\hat s)\sigma^2/n_t$, about 0.024 in M2 and M4 ($\sigma^2\approx3$); in M3, whose residuals are heteroscedastic, it is 0.25--0.38; all exceed $\delta^2=0.01$. In Heat, Step~1 has rank 30--31.

The sampling diagnostics are as follows. In M4, (H) is accepted at rate 0.56--0.59, limited mainly by the about 1\% of N1 and N2 left in U: (U) refreshes, which change neither $\tau$ nor $c$, are accepted at a similar rate (0.61--0.63), and with the true structure both updates are always accepted, whereas moving 1\% of each N direction into U reduces both rates to about 0.60. In Heat, the NQU sampler needed $1.18$--$1.27\times10^5$ calls; both methods spend about $9.5\times10^4$ calls before the window (79\% of the NQU total), and within the window the worst ESS per call of the NQU sampler is 6.5 (posterior axes) and 5.8 (original coordinates) times that of DILI. The acceptance rates of the NQU sampler there are 0.24 (Q), 0.95 (U) and 0.95 (H); simulating the slice chain on $p(\eta\mid y)$ alone gives $\IACT_{\mathcal K_H}=1.217$ for $\tau$ and 1.663 for $(\tau-\E\tau)^2$; and the LIS of DILI has dimension 17. The Gaussian model of AS-MCMC from Stein's identity is inaccurate in the stiffest directions of Heat: for seed 100, its Hessian on the 16-dimensional active subspace has a relative error of 1.1\%, with entries off by up to 39 against curvatures up to 7699. Its (H) update, which draws the Q coordinates afresh from this model, was therefore never accepted, whereas QPHD, whose least-squares fit of the quadratic $\Phi$ is exact, accepts it at rate 0.95, as the NQU sampler does. In M3, a possible mechanism for the slow mixing of Q7, not verified by sampling, is that in the true coordinates the population slope $B$ of the Q model vanishes for Q7--Q10, whereas their posterior conditional means move with $z_1$. With $\kappa_j(z_1)=m_j+b_j(z_1^2-c_0)$, a priori $\E[\Phi q_j\mid z_1]=-\lambda_j\kappa_j(z_1)$ is even in $z_1$, so the $(q_j,z_1)$ entry $\E[\Phi q_jz_1]$ of $\Hb$ vanishes; given $z_1$, the posterior of $q_j$ is $\Nor(\lambda_j\kappa_j(z_1)/(1+\lambda_j),1/(1+\lambda_j))$, whose mean has slope $2\lambda_jb_jz_1/(1+\lambda_j)$. For Q7 ($\lambda_j=10$, $b_j=0.2$) at the posterior mean $z_1\approx1.49$ (sd 0.42) this is 0.54, or 1.8 conditional standard deviations per unit of $z_1$, so Q7 should follow N1 but the Q model does not shift it. In M3 the bottleneck of AS-MCMC on the posterior axes is the second moment of the spike N2. Apart from QPHD in Heat, the smallest ratio among the gradient-free baselines, 10.8, is that of QPHD in M4 on the posterior axes. SAFES and SAFES-P shrank their steps during burn-in, and in 31 of 32 chains in M1 and M2 their most biased quantity is the second moment of N1, whose two modes the particles did not both visit; QPHD selected $k=111$--166 dimensions because of the estimation error of its 20\,301 coefficients; and the autoregression coefficient of gpCN rose to 0.996--0.998. Diagonal NUTS crossed between the modes of N1 slowly and had divergent transitions in 4 of 8 chains in M4. AS-MCMC found the informed subspace ($k=12$ in M2--M4; $k=41$--42 in M1, leaving 1.15--1.19 nats outside).

In the methane example, the ignition delay is the time at which the temperature exceeds its initial value by 400\,K, and the peak mole fractions are taken up to twice the ignition delay. The integrator uses relative and absolute tolerances $10^{-6}$ and $10^{-12}$. If a mixture does not ignite within 1\,s or the integration fails, a fixed truncation rule (ignition delay 1\,s, peaks $10^{-12}$) keeps $\Phi$ finite; this never happened in the runs reported here. One chain of $5\times10^5$ calls takes about 2.3 days on one CPU core. The Laplace approximation uses the mode found by Gauss--Newton with a Jacobian by central differences; its curvatures above 1 are 29, 10.7, 5.4, 1.45 and 1.29, and along the direction of the largest curvature the posterior has excess kurtosis about 2. The second moments use the test function $(y-c)^2$ with $c$ a fixed approximation of the posterior mean, the posterior axes are ordered by increasing variance, and the standard errors in $z$ come from the respective ESS. In the NQU chains, Step~1 has rank 3, the design of Step~2 is centred at $\mu_\cS$ with $|\mu_\cS|=4.62$ and 4.75 for the two seeds, and the acceptance rates are 0.23 (Q) and 0.21 (U). Full-space RAM starts from a prior draw with proposal covariance $2.38^2I/d$ and target acceptance 0.234, adapts throughout and discards the first 10\% ($5\times10^4$ calls) as burn-in; it was accepted at rate 0.238--0.239, and its proposal standard deviations adapted to 0.0074--0.38. Its two chains pooled largely agree with the NQU sampler, with median ratio of second moments 0.99 and ratios 0.96--1.15 on the 10 most informative posterior axes, but between its two chains $\hat R$ reaches 1.16 (posterior axes) and 1.19 (original coordinates), with 113 and 116 coordinates above 1.01. For SAFES, during burn-in $\beta$ fell to $3.4\times10^{-24}$ and $2.1\times10^{-17}$, with acceptance 0.137 and 0.146, below the lower end 0.15 of the target interval throughout. As $\beta\to0$ only the ensemble move remains and SAFES reduces to the affine-invariant sampler of Section~3.2 of \citet{dunlop2022}, whose particles cannot leave the affine hull of the ensemble (here of dimension 39) at the end of burn-in; in the sampling phase the share of the particles' spread outside that hull was about $5\times10^{-30}$. Taken at face value, the ESS of SAFES in Table~\ref{tab:cantera} implies about 13--15 times the NQU calls. Under the tuning rule of the paper ($\lambda$ fixed, $\beta$ tuned to an acceptance rate in $(0.15,0.3)$), whenever the ensemble move alone is accepted at a rate below 0.15, any controller that follows this rule drives $\beta$ to 0; this also happened in M1--M4 and Heat ($10^{-31}$--$10^{-24}$).

\end{appendices}

\makeatletter\immediate\write\@auxout{\string\citation{settings}}\makeatother
\bibliography{refs}

\end{document}